\documentclass[12pt,a4paper]{article}
\usepackage[utf8]{inputenc}
\usepackage[T1]{fontenc}
\usepackage{amsmath}
\usepackage{amsfonts}
\usepackage{amssymb}
\usepackage{graphicx}
\usepackage{parskip}
\usepackage{bbm}
\author{To define}
\usepackage{amsmath,amsfonts,amsthm,amssymb,amscd, appendix, verbatim}
\usepackage{xpatch}
\usepackage{mathtools}
\mathtoolsset{showonlyrefs}
\usepackage{setspace, color}
\usepackage{array}
\usepackage{cite}
\usepackage{hyperref}

\newtheorem{theorem}{Theorem}[section]
\newtheorem{corollary}[theorem]{Corollary}

\newtheorem{lemma}[theorem]{Lemma}

\theoremstyle{definition}

\newtheorem{remark}[theorem]{Remark}

\makeatletter
\def\@remark{\textit}
\makeatother

\newcommand{\dem}{\noindent \textit{Proof:  }}
\newcommand{\1}{\mathbbm{1}}

\newcommand{\Z}{\mathbb{Z}}
\newcommand{\F}{\mathcal{F}}
\newcommand{\Sc}{\mathcal{S}}

\newcommand{\M}{\mathcal{M}}
\newcommand{\N}{\mathbb{N}}

\newcommand{\R}{\mathbb{R}}

\newcommand{\supp}{\operatorname{supp}}

\makeatletter
\xpatchcmd{\@thm}{\fontseries\mddefault\upshape}{}{}{} % same font as thm-header
\makeatother

\usepackage[margin=1in]{geometry}

	\title{The Cauchy problem for nonlinear dispersive models of long internal waves in the presence of the Coriolis force}
 
\author{Ricardo Freire\footnote{Departamento de Ingeniería Matemática, Facultad de Ciencias Físicas y Matemáticas, Universidad de Chile, Santiago, Chile, e-mail: ricardofreiremath@gmail.com. This author was partially supported by ANID project FONDECYT 3230256.} \,
Thyago S.R. Santos \footnote{Departamento de Matemática,  Instituto de Matemática, Estatística e Computação Científica (IMECC), Universidade Estadual de Campinas (UNICAMP), Rua Sérgio Buarque de Holanda 651, Campinas, SP 13083-859, Brazil e-mail: thysou@gmail. This author was partially supported by CNPq grant 140993/2021-5} \,}
\date{}

\begin{document}
\maketitle
\begin{abstract}
\noindent We investigate models of dispersive long internal waves with rotational effects, specifically the Benjamin-Ono (BO) and intermediate long wave (ILW) equations modified by the presence of the nonlocal operator $\partial_x^{-1}$, which mathematically accounts for rotational influences. We establish a local and global well-posedness theory while ensuring the unconditional uniqueness of solutions in low-regularity Sobolev-type spaces. Our approach is based on techniques introduced by Molinet and Vento in \cite{MR3397003}.  

\end{abstract}

	
	%\tableofcontents

	\section{Introduction}
The study of nonlinear nonlocal dispersive equations has garnered significant interest in recent years due to their relevance in understanding some physical phenomena. In this work, we begin our analysis with a particular model proposed by R. Grimshaw \cite{MR797556} for describing the propagation of internal waves in a rotating fluid medium. The model is represented by the following equation, commonly referred to in the literature as the \textit{Rotational Modified Benjamin-Ono equation}
	\begin{equation}\label{BO}
		\partial_x(\partial_t u + \beta \mathcal{H} \partial_x^2 u + u \partial_xu)=  \gamma u, \quad (t,x)  \in \R \times \R, \\
    \end{equation}
% {and}
% \begin{equation}\label{ILW}
% 	\partial_x(\partial_t u -\beta \mathcal{T}_\delta \partial_x^2 u+  \partial_x u +  u \partial_xu):=  \gamma u \\
%\end{equation}
with $u=u(t,x) \in \R$. $\mathcal{H}$ denotes the Hilbert transform
\begin{equation}\label{hilberttranform}
(\mathcal{H} f) (x):= \left(-i \, \text{sgn}(\xi) \widehat{f}(\xi)\right)^\vee(x)
\end{equation}
$\beta\in \R\setminus \{0\}$ determines the type of dispersion and $\gamma \geq 0$ measures the effect of rotation in the fluid, in other words, the Coriolis force. 

The Coriolis force is a physical phenomenon that influences the motion of objects within a rotating reference frame relative to an inertial one. Commonly referred to as the \lq \lq Coriolis effect",  it describes the apparent deflection of an object's trajectory due to this rotational influence. This effect plays a crucial role in the study of internal waves in rotating fluids. Nonlinear dispersive equations incorporating Coriolis terms are frequently used in mathematical modeling to capture the interplay of dispersion, nonlinearity, and rotation, typically while neglecting dissipation. See \cite{MR797556},  \cite{MR2948883}, \cite{grimshaw1998long}  and references therein.

To facilitate the analysis, it is helpful to rewrite equation \eqref{BO} in terms of the antiderivative operator. To accomplish this, let us consider the antiderivative of a function $f$, defined as follows:
\begin{equation}\label{antiderivative}
\partial_x^{-1}f(x):=\left(-i \xi^{-1}\widehat{f}(\xi)\right)^\vee(x).
\end{equation}
Plugging \eqref{antiderivative} into \eqref{BO} we obtain the following equation
 \begin{equation}\label{BO4}
		\partial_t u + \beta \mathcal{H} \partial_x^2 u + u \partial_xu=  \gamma \partial_x^{-1} u.
    \end{equation}   
In other words, \eqref{BO} can be interpreted as the Benjamin-Ono equation with an extra nonlocal term $\partial_x^{-1} u$ added. Assuming, without loss of generality, that $\beta= 1$ we are concerned with the following IVP
\begin{equation}\label{BO2}
	\left\{
	\begin{array}{l}
	\partial_t u + \mathcal{H} \partial_x^2 u + u\partial_xu= \gamma  \partial_x^{-1} u, \quad (t,x)  \in \R \times \R,\\
		u(0,x)=u_0(x).  \tag{RMBO}
	\end{array}
	\right.
\end{equation}
Regarding the well-posedness problem for the Benjamin-Ono equation, \textit{i.e.}, when $\gamma= 0$ in \eqref{BO2}. Molinet, Saut, and Tzevtkov  in  \cite{MR1885293} showed that there cannot exist a space $X_T$ continuously embedded in $C_T H^s_x$, for any $s\in \R$, from which it is possible to obtain the necessary estimates to implement an iterative Picard scheme for the associated integral equation. See \cite{MR2172940}, \cite{MR1885293} and references therein. 

To overturn this difficulty, several compactness methods and ideas have been employed to establish a well-posedness theory for the case $\gamma=0$. For instance, Iorio \cite{MR847994} established local well-posedness in $H^s_x(\R)$ for $s> \frac{3}{2}$ using a priori energy estimates, while Ponce \cite{MR1097916} extended the result to $s=\frac{3}{2}$. Koch and Tzvetkov \cite{MR1976047} improved the local well-posedness result for $s> \frac{5}{4}$, and Kenig and Koenig \cite{MR2025062} showed local well-posedness for $s>\frac{9}{8}$, improving the results of \cite{MR1976047}. Tao, in \cite{MR2052470}, established global well-posedness for $s \geq 1$ using a gauge transformation which reduces the bad interactions in the equation. Burq and Planchon \cite{MR2357995} and Ionescu and Kenig \cite{MR2291918} proved the well-posedness in $H^{\frac{1}{4}+}_x(\R)$ and $L^2_x(\R)$  respectively, by using the gauge transformation in spaces with low regularity (see also \cite{MR2970711} for more details).

More recently, in \cite{MR4743514}, Killip, Laurens and Visan established the sharp local well-posedness of the Cauchy problem in $ H^s_x $ for $ s > -\frac{1}{2} $, in the case $ \gamma = 0 $, where the complete integrability of the BO equation is essential in their argument.

For the case $\gamma>0$, it is natural to consider the following appropriate closed subspace of $H^s_x(\R)$ in order to control the antiderivative term 
\begin{equation}\label{Espaço Zs}
Z^s_x(\R):= \left\{f \in H^s_x(\R): \|f\|_{Z^s_x}:= \|f\|_{H^s_x} + \|\partial_x^{-1} f \|_{L^2_x} <\infty\right\}.
\end{equation}
In the same spirit as the BO equation, if we consider $\gamma > 0$ in \eqref{BO2}, it is not possible to obtain the necessary estimates to implement an iterative Picard scheme in $Z^s_x(\R)$ for the associated integral equation. More precisely, based in \cite{MOSAUTTZVEKP1}, our first result ensures that there does not exist a $T > 0$ such that the problem \eqref{BO2}  satisfy the following conditions: (1) admit a unique local solution defined on the interval $[-T, T]$, and (2) the data-solution flow map $u_0 \mapsto u(t)$, $t \in [-T, T]$, is $C^2$-differentiable at the origin from $Z^s_x(\mathbb{R})$ to $Z^s_x(\mathbb{R})$.

\begin{theorem}\label{THM1 Ill-posedness}
Let $s \in \R$ and $T>0$ . Then there does not exist a space $X_T$ continuously embedded in $C([-T, T], Z^s_x(\R))$ such that
\begin{equation}\label{aIn.Ill1}
\|V(t)u_0\|_{X_T} \lesssim  \|u_0\|_{Z_x^s}, \quad u_0 \in Z_x^s(\mathbb{R}), 
\end{equation}
and
\begin{equation}\label{aIn.Ill2}
\left\| \int_0^t V(t - t') \big[u(t') \partial_x u(t')\big] \, dt' \right\|_{X_T} \lesssim  \|u\|_{X_T}^2, \quad u \in X_T. 
\end{equation}
where $\{V(t)\}_{t \in \R}$ the free evolution group associated to equation \eqref{BO2}.
\end{theorem}
As a direct consequence of the previous result, we obtain the following outcome.
\begin{corollary}
Fix $s \in \mathbb{R}$. Then there does not exist a $T > 0$ such that \eqref{BO2} admits a unique local solution defined on the interval $[-T, T]$ and such that the flow-map data-solution
$$
u_0 \mapsto u(t), \quad t \in [-T, T],
$$
for \eqref{BO2} is $C^2$-differentiable at zero from $Z^s_x(\mathbb{R})$ to $Z^s_x(\mathbb{R})$.
\end{corollary}

As far as we know, the equation in \eqref{BO2} is not completely integrable. We neither are aware of the existence of a gauge transformation as that for the BO equation ($\gamma =0$) in our case. Consequently, energy methods become essential to establish the existence of solutions for the IVP \eqref{BO2}. Linares and Milanés in \cite{MR2130214} established several properties of the solutions of the IVP \eqref{BO2}  when $\gamma>0$. In particular, a local well-posedness result for regular data can be derived using parabolic regularization techniques.

% For a complete list of references concerning local and global well-posedness for the BO equation see \cite{MR2357995, MR847994,MR2025062, MR1976047, MR1097916, MR2052470, MR2291918, MR3948114,MR2970711} and their references for further details.

By the previous remarks we are not allow at first to apply techniques utilized to study the IVP associated to the BO equation for data in $L^2_x(\R)$ or larger Sobolev spaces. However, Molinet and Vento in \cite{MR3397003} introduced a systematic method to study local well-posed for a family of one-dimensional nonlinear dispersive equations for data in $H^s_x(\R)$ for $s \geq 1/2$. The model investigated here does not fit the hypothesis of the theory developed in \cite{MR3397003}. Nevertheless we can follow their ideas and extend their method in our case.

This method developed in \cite{MR3397003}  is based on the classical energy estimate for the dyadic piece $u_N:=P_Nu$ localized around the spatial frequency $N \in 2^{\Z}$. This estimate is expressed as
\begin{equation}\label{energia introdução}
    \| u_N\|_{L^\infty_TH^s_x}^2 \lesssim \| u_N(0) \|_{H^s_x}^2 + \sup_{t \in (0,T)} \langle N\rangle^{2s} \left|\int_0^t \int_\R P_N(\partial_x u^2)(t) u_N(t) \, dx \, dt \right|.
\end{equation}
To handle the final term on the right-hand side, a standard paraproduct decomposition is used
$$
P_N(u^2) = P_N(u_{{\gtrsim N}} u_{\gtrsim N}) + 2u_{\ll N}u_N + 2N^{-1}\Pi(\partial_x u_{\ll N},u),
$$
where $\Pi$ is a pseudoproduct operator (see definition in \eqref{operator Pi}) and the derivative in \eqref{energia introdução} is moved to the lowest spatial frequencies using commutator estimates, given the frequency localized functions. This approach involves a dyadic decomposition of each function based on its modulation variable, with the one exhibiting the highest modulation placed into the Bourgain space $X^{s-1,1}$. However, addressing the characteristic function $\1_{(0,T)}$ becomes a challenge after extending the functions to $\R^2$, as it lacks continuity in $X^{s-1,1}$. Nevertheless, controlling the $X^{s-1,1}$ norm of $u$ remains relatively straightforward using classical linear estimates in Bourgain’s spaces, since for $s>1/2$ the bilinear estimates become simple (see also \cite{MR3906854}).

Using this method, the authors in  \cite{MR3397003} proved that the IVP associated with a class of strongly non-resonant dispersive equations is locally well-posed in $H^s_x(\R)$ for $s\geq 1/2$ without using gauge transformations. Therefore, our  main goal is to employ such techniques to enhance previous results about IVP \eqref{BO2}. The primary obstruction in this approach is the symbol of the equation in \eqref{BO2}
$$
\phi(\xi) :=  \xi |\xi| - \gamma\xi^{-1}
$$
which does not belong to the class of symbols considered in \cite{MR3397003} and the present analysis shall be more involved that the previous cases. In fact, the presence of the additional non-local operator $\partial_x^{-1}$ induces the formation of singularities in the resonance function associated with the symbol $\phi$ and this is a barrier to the direct application of the method.

Furthermore, considering the space $Z^s_x(\R)$ in our estimates, it results in an additional term in the energy estimate, namely
% In fact, a classical example of a function belonging to these spaces is provided by $f(x) = xe^{-x^2}$.
\begin{equation*}
J =  \sum_{N \in 2^\Z}\sup_{t \in (0,T)} \left|\int_0^t \int_\R P_N(u^2)(t) P_N \partial_x^{-1}u (t) \, dx \, dt \right|.
\end{equation*}
Next, we present the first well-posedness result of this paper.
\begin{theorem}[Local well-posedness for \eqref{BO2}]\label{Teorema de existencia}
Let $s> \frac{1}{2}$ and $\gamma >0$. Then the Cauchy problem \eqref{BO2} is unconditionally locally well-posed in $Z^s_x(\R)$, in other words, for any initial data $u_0 \in Z^s_x(\R)$, there exist a maximal time $T=T(\|u_0\|_{Z^s_x})>0$ and a unique solution to the IVP \eqref{BO2} such that $u \in C_TZ^s_x$. Furthermore, the data-solution map $u_0 \mapsto u \in C_TZ^s_x$ is continuous.
\end{theorem}
\begin{remark}
Once again, it is important to point out that the proof of Theorem \ref{Teorema de existencia} relies on compactness arguments without the use of gauge transformations (See \cite{MR2052470}). Similarly to the case of the Benjamin-Ono equation, perhaps the results obtained here (maybe) could be enhanced using such techniques. We refer to \cite{MR2291918, MR3948114, MR2970711} for further details.
\end{remark}
Since the complete integrability of the equation in \eqref{BO2} is not known for $\gamma>0$, we can leverage certain conserved quantities of the BO-flow (case $\gamma=0$) to establish the following result.
\begin{theorem}[Global well-posedness for \eqref{BO2}]\label{teorema de existencia global RMBO}
If we consider $u_0 \in Z^s_x(\R)$ with $s \geq 1$ in Theorem \ref{Teorema de existencia}, then the IVP \eqref{BO2} is unconditionally globally well-posed.
\end{theorem}

\begin{remark}
It is worth emphasizing that all the results presented here can be extended to the following class of dispersion perturbation of Benjamin-Ono equation with Coriolis term:  
\[
\partial_t u + D^\alpha_xu_x + u\partial_x u = \gamma \partial_x^{-1}u,
\]  
with $u = u(t,x)$ is real-valued, $\alpha \in [1,2]$, $\gamma>0$ and ${D^\alpha_x}$ represents the Riesz potential, defined through the Fourier transform as follows:
 $$
 {D^\alpha_x f}:= \left(|\xi|^\alpha \widehat{f}(\xi)\right)^\vee={(\mathcal{H}\partial_x)^{\alpha}f}.
 $$
\end{remark}

Next, we introduce the model proposed by Cullen and Ivanov in \cite{MR4438758}, which describes intermediate long wave propagation in dispersive media while accounting for rotational effects:
\begin{equation}\label{ILW}
	\partial_x(\partial_t v +\beta \mathcal{T}_\delta \partial_x^2 v+\frac{1}{\delta}\partial_x v +  v \partial_xv)=  \gamma v\\
\end{equation}
where
$$
\mathcal{T}_\delta f (x) := - \text{\textbf{p.v.}}\, \frac{1}{2\delta}\int_\R \coth\left(\frac{\pi(x-y)}{2\delta}\right)f(y) \, dy,
$$
with $\delta>0$ being a small non-dimensional parameter. Again, assuming, without loss of generality, that $\beta=  1$ and applying the antiderivative operator in equation \eqref{ILW}, we obtain the following equation
\begin{equation}\label{ILW4}
    \partial_t v + \mathcal{T}_\delta \partial_x^2 v+\frac{1}{\delta}\partial_x v+ v\partial_xv= \gamma \partial_x^{-1} v\quad (t,x)  \in \R \times \R.
\end{equation}

\begin{remark}
The operator $\mathcal{T}_\delta$, which arises in the study of water waves with finite depth, is associated with the Dirichlet-to-Neumann map for a two-dimensional, finite-depth, two-layer fluid (see \cite{MR4436142,MR2424620,chapouto2023deepwater,MR3625189}) . In the deep-water limit $(\delta \to \infty)$, the operator $\mathcal{T}_\delta$ converges to the Hilbert transform $\mathcal{H}$, corresponding to the Dirichlet-to-Neumann map in the case of infinite depth. This convergence suggests a potential connection  between \eqref{ILW4} and the equation in \eqref{BO2}. See \cite{chapouto2023deepwater} for more details about this discussion.
\end{remark}
We will be concerned with the following IVP associated to equation \eqref{ILW4}
\begin{equation}\label{ILW2}
	\left\{
	\begin{array}{l}
		\partial_t v + \mathcal{T}_\delta \partial_x^2 v+ \frac{1}{\delta}\partial_x v+ v\partial_xv=  \gamma\partial_x^{-1} v\quad (t,x)  \in \R \times \R,\\
		v(0,x)=v_0(x).\tag{RMILW}
	\end{array}
	\right.
\end{equation}
We will refer to equation in \eqref{ILW} as the \textit{Rotational Modified Intermediate Long Wave equation}. Note that setting $\gamma=0$ in equation \eqref{ILW} leads back to the classical ILW equation. A highly significant result in this scenario ($\gamma=0)$ was achieved recently by Chapouto, Li, Oh and Pilod in \cite{chapouto2023deepwater} (see also Ifrim and Saut  \cite{ifrim2023lifespan}), where the authors established local and global well-posedness theory for initial data in  $L^2_x(\R)$. For further information, see \cite{MR1044731},\cite{ifrim2023lifespan}, \cite{MR3397003},  \cite{MR3931835}, \cite{MR2887985}, and their references.

The symbol of the equation \eqref{ILW4} is given by
$$
\omega_\delta(\xi) := \xi^2\coth{(\delta\xi)}+ \frac{1}{\delta}\xi - \gamma \xi^{-1}.
$$
Due to the presence of the singular term $\xi^{-1}$, this symbol does not satisfy Hypothesis 1 considered by Molinet and Vento in \cite{MR3397003}. Hence, the challenges encountered in the case of \eqref{BO2} also manifest in this model.

Due to the close relationship between the IVPs \eqref{ILW2} and \eqref{BO2}, it is natural to apply the same methods used to prove Theorem \ref{Teorema de existencia} to obtain well-posedness results for the IVP \eqref{ILW2}. Consequently, by employing analogous techniques, we derive the following result.

\begin{theorem}[Local well-posedness for \eqref{ILW2}]\label{Teorema de existencia ILW}
Let $s> \frac{1}{2}$ and $\gamma >0$. Then, the Cauchy problem \eqref{ILW2} is unconditionally locally well-posed in $Z^s_x(\R)$. Furthermore, the data-solution map $u_0 \mapsto u \in C_TZ^s_x$ is continuous.
\end{theorem}

Since the ILW equation is completely integrable, we can follow the same approach used for \eqref{BO2} in this case to obtain the following global result:

\begin{theorem}[Global well-posedness for \eqref{ILW2}]\label{Teorema global ILW}
For $s\geq1$, as stated in Theorem \ref{Teorema de existencia ILW}, IVP \eqref{ILW2} is unconditionally globally well-posed.
\end{theorem}

It is important to emphasize that, up to the time of writing this work, we have not found any work focusing on a well-posedness theory for the IVP \eqref{ILW2}.

The paper is organized as follows: In Section \ref{seção das notação}, we introduce essential notations, define the necessary functional spaces, and present key tools for our arguments. Section \ref{res preliminares} contains some preliminary lemmas. In Section \ref{ill posedness}, we establish the Ill-Posedness results. Section \ref{RMBO} is dedicated to studying equation \eqref{BO2}, presenting \textit{a priori} and energy estimates, as well as proving Theorems \ref{tempo existencia} and \ref{teorema de existencia global RMBO}. Finally, Section \ref{seçao ILW} provides the proofs of our results for equation \eqref{ILW2}.

\section{Notations and basic concepts}\label{seção das notação}
For any two non-negative quantities $X$ and $Y$, the notation $X\lesssim Y$ implies the existence of an absolute constant $C>0$ such that $X \leq CY$. We will denote $X\sim Y$ when both $X \lesssim Y$ and $Y \lesssim X$. The notation $X \ll Y$ is used when $X \gtrsim Y$ does not hold. Additionally, if $\alpha \in \R$, then $\alpha_+$ and $\alpha_-$ represent numbers slightly greater and smaller than $\alpha$, respectively.

For $u \in \Sc$, where $\Sc=\Sc(\R)$ denotes the Schwartz class, and for $s \in \R$ we will define the \textit{Bessel potential of order $-s$} in variable $x$ by the multiplier
$$
\widehat{J_x^s u}(\xi) := \langle \xi \rangle^s \widehat{u}(\xi),
$$
where
$$
\widehat{f}(\xi)=\F_x(f)(\xi):= \int_\R e^{-i x \cdot \xi} f(x)\, dx,
$$
denotes the usual Fourier transform in variable and $\langle x \rangle : = (1 + |x|^2)^{\frac{1}{2}}$. We will consider $H^s_x$ the usual Sobolev space equipped with the norm
$$
\|u\|_{H^s_x}: =\|J^s_x u \|_{L^2_x}.
$$
Denote by $\{V(t)\}_{t \in \R}$ the free evolution group associated to equation \eqref{BO2}. In other words, for $f \in \Sc(\R)$
$$
\left(V(t)f\right)(x) := c \int_\R e^{i(x\cdot \xi -t \phi(\xi))} \widehat{f} (\xi) \, d\xi= \left( \exp(-it \phi(\xi)) \widehat{f}(\xi)\right)^\vee (x)
$$
where $\phi(\xi) = \xi |\xi| - \gamma\xi^{-1}$ is the symbol of the equation \eqref{BO2}.  For $s,b \in \R$ we define the Bourgain space $X^{s,b}=X^{s,b}_{\tau= \phi(\xi)}$ associated with the symbol $\phi$ as the completion of the Schwartz space $\Sc(\R \times \R)$ under the norm
$$
\|u\|_{X^{s,b}_{\tau= \phi(\xi)}}^2 := \int_\R \int_\R \langle \tau- \phi(\xi)\rangle^{2b} \langle \xi \rangle^{2s} |\tilde{u}(\tau, \xi)|^2 \, d\tau d\xi,
$$
where $\tilde{u}=\tilde{u}(\tau,\xi)$ represents the space-time Fourier transform of $u=u(t,x)$. The spaces $X^{s,b}$ was introduced by Jean Bourgain for the study of dispersive equations (see \cite{MR1209299,MR1215780}).

The definitions for \eqref{ILW2} will follow similarly, but first, consider a change of coordinates given by  
$$
\bar{v}(t,x) := \delta v(\delta^2 t, \delta x).
$$  
Thus, we can, without loss of generality, set $ \delta = 1 $ in \eqref{ILW2}. Therefore, we consider the symbol of equation \eqref{ILW2}
$$
\omega(\xi):=\omega_1(\xi)=  \xi^2\coth{(\xi)}+ \xi - \gamma\xi^{-1}
$$  
and $ \{W(t)\}_{t\in \R}$ the evolution group associated to \eqref{ILW2}, \textit{i.e.}
$$
\left(W(t)g\right)(x) := c \int_\R e^{i(x\cdot \xi -t \omega(\xi))} \widehat{g} (\xi) \, d\xi =\left( \exp(-it \omega(\xi)) \widehat{g}(\xi)\right)^\vee (x).
$$
When there is no risk of ambiguity, we will use the notation $ X^{s,b} $ to represent the Bourgain space associated simultaneously with the symbols $ \phi(\xi)$ and $ \omega(\xi)$.
% \begin{remark}
% It is important to emphasize again that the symbols $\phi(\xi)$ and $\omega(\xi)$ considered in this work do not belong to the class of symbols described in Hypothesis 1 of \cite{MR3397003}, primarily due to the presence of the singular term $\xi^{-1}$.
% \end{remark}

One of the key tools employed in this work is the Littlewood-Paley decomposition. Therefore, let us fix a smooth cutoff function $\eta \in C^\infty_c(\mathbb{R})$ satisfying the following properties: 
\begin{enumerate}
    \item $\supp \eta \subset [-2,2]$;
    \item $0 \leq \eta \leq 1$ on $[-2,2]$;
    \item $\eta \equiv 1$ on $[-1,1]$.
\end{enumerate}
Define $\varphi(\xi) := \eta(\xi) - \eta(2\xi)$, and for a dyadic number $N \in 2^\mathbb{Z} := \{2^k : k \in \mathbb{Z}\}$, let $\varphi_N$ be given by  
$$
\varphi_N(\xi) := \varphi(N^{-1} \xi).  
$$  
With these definition we have the following relation 
\begin{equation*}
\supp \varphi_N \subset \{\xi \in \R: |\xi| \sim N\}, \quad\text{with} \quad\sum_{N\, \in\, 2^{\Z}} \varphi_N(\xi) = 1, \quad \forall \,\, \xi \in \R^{*}.
\end{equation*}
For each dyadic number $N \in 2^\mathbb{Z}$ and  $f \in \mathcal{S}(\mathbb{R})$ we define the \textit{Littlewood-Paley multipliers} as follows:    
$$
\widehat{P_N f}(\xi) :=  \varphi_N(\xi) \widehat{f}(\xi), \,\,\widehat{P_{>N} f}(\xi) :=  \left(1-\varphi_N(\xi)\right) \widehat{f}(\xi) \,\, \text{and}\,\,  \widehat{P_{\leq N} f}(\xi) :=  \eta(N^{-1}\xi) \widehat{f}(\xi).
$$    
Using these definitions, we have the following summation identities 
$$
P_{\geq N} f(x) := \sum_{M \geq N} P_M f(x), \quad P_{\leq N} f(x) := \sum_{M \leq N} P_M f(x),
$$  
where $M$ ranges in $2^\mathbb{Z}$. Using the summation identities defined at the beginning of this section, we can define the operators $P_{\ll N}$, $P_{\gg N}$, $P_{\lesssim N}$, and others. We now present the Littlewood-Paley theorem, which provides a estimate that is extensively utilized in this work.

\begin{theorem}[See \cite{Muscalu_Schlag_2013}]\label{Littlewood-Paley Theorem}
Let $1<p\leq \infty$. Then for $f \in \Sc(\R)$
$$
 \| f\|_{L^p_x} \sim \| (\sum_{N \in 2^{\Z}} |P_N f|^2)^{1/2} \|_{L^p_x}.
$$
\end{theorem}
Considering $p=2$ and by  the Plancherel identity we have the identity
$$
\|f\|_{H^s_x}^2 \sim \sum_{N \in 2^{\Z}} \langle N\rangle^{2s} \|P_N f\|_{L^2_x}^2 \sim \|P_{\leq 1} f\|_{L^2_x}^2 + \sum_{N >1} N^{2s} \|P_N f\|_{L^2_x}^2.
$$
In conclusion, let us cite the Lemma concerning the nonlinear estimates in Sobolev spaces.
\begin{lemma}[See \cite{MR2424078}]\label{lemma do sobolev}
Let $(s_1,s_2,s_3) \in \R^3$ such that
$$
s_1 \geq s_3, \,\, s_2 \geq s_3, \,\,\, s_1 + s_2 \geq 0 \,\,\, \text{and}\,\,\, s_1 + s_2 -s_3 > \frac{1}{2}
$$
and $f_i \in H^{s_i}_x(\R)$ for $i=1,2$. Then $f_1 f_2 \in H^{s_3}_x(\R)$ and
\begin{equation}
\|f_1f_2\|_{H^{s_3}_x} \lesssim \|f_1\|_{H^{s_1}_x}\|f_2\|_{H^{s_2}_x}.
\end{equation}  
\end{lemma}

\section{Preliminary results} \label{res preliminares}
Throughout this work, we will use time-restricted versions of certain functional spaces. For $ T > 0 $, we define the restricted space $ B_T $, where $ B $ is a functional space, as the space of space-time functions $ f: (0,T) \times \mathbb{R} \to \mathbb{R} $ that satisfy the following condition:
$$
\| f \|_{B_T} := \inf \left\{\|F\|_{B} \,\, | \,\,  F: \R \times \R \rightarrow \R,\,\, F_{\mid_{ (0,T)\times \R}} \equiv f \right\} < \infty.
$$
From this point on, it will also be necessary to consider a space-time version of the Littlewood-Paley multipliers, which localize temporally on the hypersurfaces $ \{\tau = \phi(\xi)\} $ and $ \{\tau = \bar{\phi}(\xi)\} $. Thus, for $ L \in 2^{\mathbb{N}} := \{2^n : n \in \mathbb{N}\} $, we define
$$
\psi_L(\tau, \xi):=\varphi \left(L^{-1}(\tau- \phi(\xi))\right)
$$
and 
$$ \bar{\psi}_L(\tau, \xi)=\varphi \left(L^{-1}(\tau- \bar{\phi}(\xi))\right), $$
where $\varphi$ was defined in the Section \ref{seção das notação}. With these definitions we have the following relation 
\begin{equation*}
\supp \psi_L \subset \{(\tau,\xi) \in \R^2 : | \tau- \phi(\xi) | \sim L\}, \quad\text{with}\quad \sum_{L \,\in \, 2^{\N}} \psi_L (\tau, \xi) = 1, \quad \forall \,\, (\tau, \xi) \in \R^2 .
\end{equation*}
 For each dyadic number $L \in 2^\mathbb{N}$, we define the {Littlewood-Paley multipliers} as follows:  
$$
Q_L g(t,x) := \mathcal{F}_{t,x}^{-1} \left( \psi_L(\tau, \xi) \tilde{g}(\tau, \xi) \right)(t,x),
$$  
$$
\bar{Q}_L g(t,x) := \mathcal{F}_{t,x}^{-1} \left( \bar{\psi}_L(\tau, \xi) \tilde{g}(\tau, \xi) \right)(t,x)
$$  
where $g \in \Sc(\mathbb{R}\times \R)$. As done previously, it is naturally possible to define the operators $ Q_{\leq L}, Q_{\ll L}, \bar{Q}_{\geq L} $, and so on.
\begin{lemma}[See \cite{MR3397003}]
Let $L \in 2^{\N}$, $1 \leq p \leq \infty$ and $s \in \R$. The operators $Q_{\leq L}$ and $\bar{Q}_{\leq L}$ are uniformly continuous in the space $L^p_t H^s_x$.\qed
\end{lemma}
It will be crucial to extend functions defined on the interval $[0,T]$ to the whole real line $\R$ in order to work with Bourgain spaces. To this end, we need to handle the characteristic function $\1_T$ of the interval $[0,T]$, which is generally not continuous in such spaces. Therefore, for a given $R > 0$, let us consider the decomposition:
$$
\1_T= \1_{T,R}^{\text{low}} + \1_{T,R}^{\text{high}},
$$
where 
$$
\F_t(\1_{T,R}^{\text{low}})(\tau)= \eta(\frac{\tau}{R})\F_t(\1_T)(\tau)\quad \text{and}\quad \F_t(\1_{T,R}^{\text{high}})(\tau)=(1- \eta(\frac{\tau}{R}))\F_t(\1_T)(\tau).
$$
The properties of this decomposition are provided in the following lemma.
\begin{lemma}[See \cite{MR3397003, palacios2022local}] \label{Lemma 2.4}
For any $R>0$ and  $T>0$, we have
\begin{equation}\label{2.9}
	\|\1_{T,R}^{\text{low}}\|_{L^\infty_t} \lesssim 1,
\end{equation}
and,  for $1 \leq q \leq \infty,$
\begin{equation}\label{2.8}
	\|\1_{T,R}^{\text{high}}\|_{L^q_t} \lesssim \min\left\{T, \frac{1}{R}\right\}^{1/q}.
\end{equation}
Moreover, if $f \in L^2_{t,x}(\R^2)$ and $L \in 2^{\N}$ such that $L \gg R$ we have
$$
\|Q_L(\1_{T,R}^{\text{low}}f) \|_{L^2_{t,x}} \lesssim \|Q_{\sim L}f \|_{L^2_{t,x}}.
$$
and
$$
\|\bar{Q}_L(\1_{T,R}^{\text{low}}f) \|_{L^2_{t,x}} \lesssim \|\bar{Q}_{\sim L}f \|_{L^2_{t,x}}.
$$\qed
\end{lemma}

\section{Ill-posedness result}\label{ill posedness}

This subsection is dedicated to the proof of Theorem \ref{THM1 Ill-posedness}. For the sake of completeness, we will emphasize the key differences between this proof and the result presented by Molinet, Saut, and Tzvetkov in \cite{MR1885293}. Our focus will be on the case of \eqref{BO2}. For the other case, the equation \eqref{ILW4}, we suggest following the same steps detailed here and consulting the ideas discussed in Section 3 of \cite{MR1885293}.

The proof will proceed by contradiction. First, we choose a function $ u_0 $ designed to contradict the assumptions of the theorem. More precisely, $ u_0 $ is defined through its Fourier transform as follows
$$
\widehat{u_0}(\xi) = \alpha^{-\frac{1}{2}} \1_{I_1}(\xi) + \alpha^{-\frac{1}{2}} N^{-s} \1_{I_2}(\xi),
$$
where $I_{\alpha} = \left[\frac{\alpha}{2}, \alpha \right]$, $I_N = [N, N + \alpha]$, $N \gg 1$ and $0 < \alpha \ll 1$. Note that $u_0 \in Z_x^s(\R)$ and 
$$
\|u_0\|_{Z_x^s} \sim  \frac{3}{2} + \frac{1}{\alpha} + \frac{1}{N^{\frac{2s+1}{2}}}.
$$ 
Now, using the symbol $\phi$ of \eqref{B04} and the inverse Fourier transform, we obtain 
\begin{equation}
    \begin{aligned}
        \int_0^t V(t - t')& \big( (V(t')u_0)(V(t')\partial_x u_0) \big) \, dt'=\\ 
&= c \int_0^t \int_{\mathbb{R}} e^{ix\xi + it\phi(\xi)} e^{-it'\phi(\xi)\xi} 
\big( e^{it'\phi(\cdot)} \widehat{u_0}(\cdot) \big) \ast \big( e^{it'\phi(\cdot)} \widehat{u_0}(\cdot) \big)(\xi) \, d\xi \, dt'\\
&= c \int_{\mathbb{R}^2} e^{ix\xi + it\phi(\xi)} 
\frac{e^{it(\phi(\xi_1) + \phi(\xi - \xi_1) - \phi(\xi))} - 1}{\phi(\xi_1) + \phi(\xi - \xi_1) - \phi(\xi)} 
\xi \, \widehat{u_0}(\xi_1) \widehat{u_0}(\xi - \xi_1) \, d\xi_1 \, d\xi\\
& = \frac{c}{\alpha} \int \limits_{\substack{\xi_1 \in I_{\alpha}\\ \xi-\xi_1 \in I_{N} }} e^{ix\xi + it\phi(\xi)} 
\frac{e^{it(\phi(\xi_1) + \phi(\xi - \xi_1) - \phi(\xi))} - 1}{\phi(\xi_1) + \phi(\xi - \xi_1) - \phi(\xi)} 
\xi  \, d\xi_1 \, d\xi\\
&\quad +  \frac{c}{\alpha N^{2s}} \int\limits_{\substack{\xi_1 \in I_{\alpha}\\ \xi-\xi_1 \in I_{N} }} e^{ix\xi + it\phi(\xi)} 
\frac{e^{it(\phi(\xi_1) + \phi(\xi - \xi_1) - \phi(\xi))} - 1}{\phi(\xi_1) + \phi(\xi - \xi_1) - \phi(\xi)} 
\xi  \, d\xi_1 \, d\xi\\
&\quad +  \frac{c}{\alpha N^{s}} \int\limits_{\substack{\xi_1 \in I_{\alpha}\\ \xi-\xi_1 \in I_{N} }} e^{ix\xi + it\phi(\xi)} 
\frac{e^{it(\phi(\xi_1) + \phi(\xi - \xi_1) - \phi(\xi))} - 1}{\phi(\xi_1) + \phi(\xi - \xi_1) - \phi(\xi)} 
\xi  \, d\xi_1 \, d\xi\\
&\quad +  \frac{c}{\alpha N^{s}} \int \limits_{\substack{\xi_1 \in I_{\alpha}\\ \xi-\xi_1 \in I_{N} }} e^{ix\xi + it\phi(\xi)} 
\frac{e^{it(\phi(\xi_1) + \phi(\xi - \xi_1) - \phi(\xi))} - 1}{\phi(\xi_1) + \phi(\xi - \xi_1) - \phi(\xi)} 
\xi  \, d\xi_1 \, d\xi \\
&=: f_1(t,x) + f_2(t,x) + f_3(t,x) + f_4(t,x).
    \end{aligned}
\end{equation}

Assume the existence of a function space $ X_T $ such that inequalities \eqref{aIn.Ill1} and \eqref{aIn.Ill2} are satisfied. Substituting $ u = V(t)u_0 $ into inequality (8), it follows that
$$
\left\|\int_0^t V(t - t') \big[(V(t')u_0)(V(t')\partial_x u_0)\big] \, dt'\right\|_{X_T} \lesssim  \|V(t)u_0\|_{X_T}^2.
$$
Furthermore, by leveraging inequality \eqref{aIn.Ill2} and the fact that $ X_T $ is continuously embedded in $ C([-T, T], Z^s_x(\mathbb{R})) $, we deduce that for any $ t \in [-T, T] $,
\begin{equation}\label{ineqint}
\left\|\int_0^t V(t - t') \big[(V(t')u_0)(V(t')\partial_x u_0)\big] \, dt'\right\|_{Z_x^s} \lesssim \|u_0\|_{Z_x^s}^2. \tag{9}
\end{equation}
Now, observe that since the supports of $ \widehat{f_1}(\xi),\, \widehat{f_2}(\xi)$ and $ \widehat{f_3 + f_4}(\xi)  $ are disjoint (see \cite{MR1885293} for further details), it follows that
\begin{equation}\label{finalineq}
\left\|\int_0^t V(t - t')[ (V(t')u_0)(V(t')\partial_x u_0) ] \, dt'\right\|_{Z^s_x} \geq \| (f_3 +f_4)(t, \cdot) \|_{Z_x^s}. 
\end{equation}
Thus, defining 
$$
\chi(\xi,\xi_1) = \phi(\xi_1) + \phi(\xi - \xi_1) - \phi(\xi)
$$
and  using the Taylor expansion, we obtain:
$$
\left|\frac{e^{it\chi(\xi,\xi_1)}-1}{\chi(\xi,\xi_1)}\right| \sim  |t|\sqrt{1 + \frac{(t\chi(\xi,\xi_1))^2}{4}} \sim |t| + \left(|t|^2 \left( \frac{1}{\alpha} + \frac{1}{N} - 2\alpha N\right)\right)
$$
Now, after a straightforward computation, we obtain the following inequalities.
$$
\| (f_3 +f_4)(t, \cdot) \|_{H_x^s} \gtrsim \alpha^{1/2} N\cdot \left(|t|^2 \left( \frac{1}{\alpha} + \frac{1}{N} - 2\alpha N\right)\right)  
$$
and
$$
\|\partial_x^{-1} (f_3 +f_4)(t, \cdot) \|_{L_x^2} \gtrsim \frac{\alpha^{1/2}}{N^s}\cdot \left(|t|^2 \left( \frac{1}{\alpha} + \frac{1}{N} - 2\alpha N\right)\right)  
$$
Then, 
$$
\| (f_3 +f_4)(t, \cdot) \|_{Z_x^s} \gtrsim (\alpha^{1/2}N + \frac{\alpha^{1/2}}{N^s} ) \cdot \left(|t|^2 \left( \frac{1}{\alpha} + \frac{1}{N} - 2\alpha N\right)\right)  
$$
Therefore, by the \eqref{finalineq} and \eqref{ineqint} we conclude
$$
\frac{3}{2} + \frac{1}{\alpha} + \frac{1}{N^{\frac{2s+1}{2}}} \sim \|u_0\|_{Z^s} \gtrsim \| (f_3 +f_4)(t, \cdot) \|_{Z_x^s} \gtrsim (\alpha^{1/2}N + \frac{\alpha^{1/2}}{N^s} )\cdot |t|^2 \left( \frac{1}{\alpha} + \frac{1}{N} - 2\alpha N\right) 
$$
By appropriately choosing $ N \gg 1 $ and $ \alpha \ll 1 $, the proof is completed via contradiction, thus establishing the theorem.
\qed

\section{The Rotational Modified Benjamin-Ono Equation}\label{RMBO}

We will start our analysis by examining equation \eqref{BO2}, i.e., the Rotational Modified Benjamin-Ono equation. Therefore, consider the \textit{resonance function} given by the symbol $\phi(\xi)$ 
\begin{equation}\label{resonant function}
\Omega(\xi_1,\xi_2):=\phi(\xi_1 + \xi_2) - \phi(\xi_1) - \phi(\xi_2),
\end{equation}
for $(\xi_1,\xi_2) \in \R^2$. The next lemma gives us a control over the absolute value of \eqref{resonant function}.
\begin{lemma}\label{controle resonant}
Consider $(\xi_1,\xi_2) \in \R^2$. If $\min \{|\xi_1|, |\xi_2|\} \gg 1$ then
$$
| \Omega(\xi_1,\xi_2)| \lesssim |\xi|_{\max}^2|\xi|_{\min} +|\xi|_{\max}|\xi|_{\min}^{-2},
$$
where $|\xi|_{\max}= \max \{|\xi_1|,|\xi_2|, |\xi_1 + \xi_2|\}$ and $|\xi|_{\min}= \min \{|\xi_1|,|\xi_2|, |\xi_1 + \xi_2|\}.$
\end{lemma}
% \begin{remark}
% Throughout this work, we shall employ the notation $\phi(\xi)$ to represent $\phi_+(\xi)$. Consequently, we will omit the subindex "$\pm$" in the proofs. To the case $\phi(\xi)=\phi_{-}(\xi)$ follows analogously. The same will apply to $\omega_{\pm}.$
% \end{remark}

\dem  Without loss of generality, let us assume $\gamma = 1$ throughout this proof. First, let us split the symbol as follows:
$$
\phi (\xi)= \phi_g (\xi) - \phi_b(\xi),
$$
where $\phi_g (\xi):= -\xi |\xi|$ and $\phi_b (\xi):= \frac{1}{\xi}$. Now, we can obtain the following decomposition for the resonance function
$$
\Omega(\xi_1,\xi_2) : = \Omega_g (\xi_1,\xi_2) - \Omega_b (\xi_1,\xi_2),  
$$
with 
\begin{equation}\label{omega g}
\Omega_g (\xi_1,\xi_2)= (\xi_1 + \xi_2) |\xi_1 + \xi_2| - \xi_1|\xi_1| -\xi_2|\xi_2|
\end{equation}
and
\begin{equation}\label{omega b}
	\Omega_b (\xi_1,\xi_2)= \frac{1}{\xi_1 +\xi_2} - \frac{1}{\xi_1} - \frac{1}{\xi_2}.
\end{equation}
Since that $|\xi_1|\geq \min \{|\xi_1|, |\xi_2|\}  \gg 1$, the following estimate is valid
\begin{equation}
|\Omega_g (\xi_1,\xi_2)| \lesssim |\xi|_{\max}^2|\xi|_{\min}.
\end{equation}
Thus, if we obtain an estimate for $|\Omega_b (\xi_1,\xi_2)|$, the result follows. By symmetry, we can assume $|\xi_2| \leq |\xi_1|$ and separate our analysis in some cases:

\noindent \textbf
{{Case 1}:}   $|\xi_2| \ll |\xi_1|$. Since $|\xi_1| \gg 1$, there exists $ \theta >0$ with $\xi_1\leq \theta \leq \xi_1 + \xi_2$ such that
\begin{equation}\label{estimativa soma}
\left| \frac{1}{\xi_1 +\xi_2} - \frac{1}{\xi_1} \right| = |\xi_2| \left|\theta\right| ^{-2}\lesssim |\xi_2| \left|\xi_1\right| ^{-2}.
\end{equation}
Now, since $|\xi_2| \leq |\xi_1|$ we have, by hypothesis, $1 \leq |\xi_1||\xi_2|^{-1}$, then
\begin{equation}\label{estimativa xi_2}
\frac{1}{|\xi_2|} \lesssim \frac{|\xi_1||\xi_2|^{-1}}{|\xi_2|} = |\xi_1||\xi_2|^{-2}.
\end{equation}
Gathering \eqref{estimativa soma} and \eqref{estimativa xi_2} leads to
\begin{equation}
|\Omega_b (\xi_1,\xi_2)| \lesssim |\xi_1||\xi_2|^{-2} + |\xi_2| \left|\xi_1\right| ^{-2}  \lesssim |\xi|_{\max}|\xi|_{\min}^{-2}. 
\end{equation}

\noindent \textbf{{Case 2}}: $|\xi_2| \gtrsim  |\xi_1| $. In this case $|\xi_2| \gtrsim  |\xi_1| \gg 1 $ and, by symmetry we can assume that $\xi_1 >0$. First, suppose that $\xi_1 \xi_2 \geq 0$, then $0<1 \ll \xi_2 < \xi_1 +\xi_2$  and
\begin{align*}
|\Omega_b (\xi_1,\xi_2)| &= \left| \frac{1}{\xi_1 +\xi_2} - \frac{1}{\xi_1} - \frac{1}{\xi_2}\right| \\
&=\left|\frac{1}{\xi_1 +\xi_2} - \frac{1}{\xi_1} - \frac{1}{\xi_2} +1 - 1 + \frac{1}{\xi_1 +1 } -\frac{1}{\xi_1 +1 }\right| \\
%&= \left|\int_1^{\xi_2} \frac{d}{dt} \left( \frac{1}{\xi_1 + t} - \frac{1}{t}\right)\, dt + \left( \frac{1}{\xi_1 + 1} - \frac{1}{\xi_1}\right) - 1\right| \\
&=\left| \int_1^{\xi_2} \left( \frac{1}{t^2} -\frac{1}{(\xi_1 + t)^2} \right)\, dt + \left( \frac{1}{\xi_1 + 1} - \frac{1}{\xi_1}\right) - 1\right| \\
%& \leq  \int_1^{\xi_2} \left| \frac{1}{t^2} -\frac{1}{(\xi_1 + t)^2} \right|\, dt + \left|  \frac{1}{\xi_1 + 1} - \frac{1}{\xi_1}\right| + 1 \\
& \lesssim  \xi_1^{-2} \int_1^{\xi_2}\, dt + \xi_2\xi_1^{-2} + 1\\
%& \lesssim \xi_1^{-2} \xi_2 - \xi_1^{-2} + \xi_2\xi_1^{-2} + 1 \\
& \lesssim  |\xi|_{\max}|\xi|_{\min}^{-2}.
\end{align*}

\noindent\textbf{{Case 3}}: $|\xi_1| \sim |\xi_2|$.  When $ \xi_1 $ is approximately equal to $ -\xi_2 $, a singularity occurs, making it essential to consider this specific case. In this scenario, $ |\xi_1 + \xi_2| $ becomes very small but is controlled by $ |\xi_i| $ for $ i = 1, 2 $. Thus, we can deduce that
$$
|\Omega_b(\xi_1,\xi_2)| \lesssim \frac{1}{|\xi_1+\xi_2|} \lesssim |\xi|_{\max}|\xi|_{\min}^{-2}.
$$
For the situation where $\xi_1 \xi_2 < 0$ we separate in the cases $\xi_1 + \xi_2 \ll -\xi_2$ and	$\xi_1 + \xi_2 \gtrsim -\xi_2$ and  we can argue exactly as in the previous cases.

\qed

Now, let $\Psi \in L_x^\infty(\mathbb{R}^2)$ and define the \textit{pseudoproduct operator} $\Pi_\Psi$ as:
\begin{equation}\label{operator Pi}
\widehat{\Pi_\Psi (f,g)}(\xi)=\int_\R \widehat{f}(\zeta)\widehat{g}(\xi - \zeta) \Psi(\xi, \zeta) \, d\zeta
\end{equation}
where $f,g \in \Sc(\R)$. Note that the expression in \eqref{operator Pi} defines a bilinear operator on $\Sc(\R) \times \Sc(\R)$ and satisfies
$$
\int_\R \Pi_\Psi(f,g)(x)h(x)\, dx=\int_\R f(x)\Pi_{\Psi_1}(g,h)(x)\, dx=\int_\R \Pi_{\Psi_2}(f,h)(x)g(x)\, dx,
$$
where $\Psi_1(\xi, \zeta):=\Psi(\zeta, \xi)$ and $\Psi_2(\xi, \zeta):=\Psi(\xi -\zeta, \xi)$, for all $f,g,h \in \Sc(\R)$. In particular, 
$$
\Pi_\Psi(f,g)(x)= f(x)g(x), \,\, \text{when} \,\, \Psi \equiv 1.
$$
Throughout this work we write $\Pi:= \Pi_\Psi$, where $\Psi$ can vary in each situation.

From the Lemma \ref{controle resonant} and the Bernstein's inequalities and observing that
$$
\Omega(\xi_1,\xi_2)=-\left[(\tau_1+\tau_2 - \phi(\xi_1 + \xi_2))- ( \tau_1- \phi(\xi_1))- (\tau_2 - \phi(\xi_2))\right]
$$
for all $\tau_1, \tau_1 \in \R$, we have, directly, a similar result as the one established in \cite{MR3397003}.
\begin{lemma}[See \cite{MR3397003}]\label{lemma 2.2}
	Let $N_i \in 2^\Z$, $L_i \in 2^\N$ with $i=1,2,3$ and $0<N_1 \leq N_2 \leq N_3$. Furthermore, consider $f_i \in \Sc(\R\times \R)$ such that
$$
\text{supp }\widehat{f_i}\subset \{\xi \in \R: |\xi| \sim N_i \},
$$
for each $i=1,2,3$. Then
$$
\int_\R\int_\R (\Pi(Q_{L_1}f_1,Q_{L_2}f_2) Q_{L_3}f_3)(t,x)\,\, dt\, dx =0,
$$
for all pseudoproduct $\Pi \equiv \Pi_{\Psi}$, whenever the $\displaystyle \max_{i=1,2,3}L_i \gtrsim N_1N_2^2 +N_1^{-2}N_2$  is not satisfied.\qed
\end{lemma}

\subsection{\textit{A Priori} Estimates}\label{a priori estimates}

First, we consider the following estimates involving trilinear pseudoproducts operators.
\begin{lemma}\label{l3.4}
Let $u_i \in L^{2}_t L ^2_x \cap X^{0,1}$, $i=1,2,3$, be functions with spatial Fourier support in the annulus $\{\xi \in \R: |\xi| \sim N_i\}$ with $N_1 \leq N_2 \leq N_3$ where $N_i \in 2^\Z$ for all $i=1,2,3.$ For any $t>0$, we set the trilinear form
$$
I_t(u_1,u_2,u_3):=\int_0^t \int_\R \Pi(u_1,u_2)(t,x)u_3(t,x)\,\, dx dt.
$$
If $N_1 \lesssim 1 :$
\begin{equation}\label{3.4}
|I_t(u_1,u_2,u_3)| \lesssim N_1^{1/2} \|u_1\|_{L^\infty_t L^2_x}\|u_2\|_{L^2_tL^2_x}\|u_3\|_{L^2_tL^2_x}.
\end{equation}
Case $N_1 \gg 1$,
\begin{align*}
|I_t(u_1,u_2,u_3)| &\lesssim N_1^{1/2} R^{\frac{1}{p}-1}\|u_1\|_{L^p_tL^2_x}  \| u_2\|_{L^\infty_tL^2_x}\| u_3\|_{L^\infty_tL^2_x}\\
&\quad + {N_1^{1/2}}L^{-1} \|{u}_1\|_{X^{0,1}} \| u_2\|_{L^2_tL^2_x}\| u_3\|_{L_t^\infty L^2_x} \\
&\quad + N_1^{1/2} L^{-1}  \|u_1\|_{L^\infty_t L^2_x} \left(\|u_2\|_{X^{0,1}}\| u_3\|_{L_t^2 L^2_x} + \| u_2\|_{L_t^2 L^2_x}  \|u_3\|_{X^{0,1}}\right),
\end{align*}
for $1 \leq p < \infty$ , $0< R \ll L$ with  $L \in 2^\N$ that does not satisfy $L \gtrsim N_1N_2^2 +N_2N_1^{-2}$.
\end{lemma}
\dem   Let $N_1 \in 2^\Z$. Suppose first  that $N_1 \lesssim 1$, then by the Bernstein's inequalities
\begin{align*}
|I_t(u_1,u_2,u_3)|&=\left|\int_0^t \int_\R \Pi(u_1,u_2)(t,x)u_3(t,x)\,\, dx dt\right| \\
&\leq \int_0^t \left|\int_\R \Pi(u_1,u_2)(t,x)u_3(t,x)\,\, dx \right| dt \\
& \leq N_1^{1/2} \int_0^t \|u_1(t)\|_{L^2_x} \|u_2(t)\|_{L^2_x} \|u_3(t)\|_{L^2_x}\, dt \\
& \leq N_1^{1/2}\|u_1\|_{L^\infty_tL^2_x}  \int_0^t \|u_2(t)\|_{L^2_x} \|u_3(t)\|_{L^2_x} \, dt\\
& \leq  N_1^{1/2} \|u_1\|_{L^\infty_tL^2_x} \|u_2\|_{L^2_tL^2_x} \|u_3\|_{L^2_tL^2_x}. 
\end{align*}
For the case $ N_1 \gg 1 $, let $ R > 0 $ and split $ I_t $ as follows:
\begin{align*}
I_t(u_1,u_2,u_3)=& \int_0^t \int_\R \Pi(u_1,u_2)(t,x)u_3(t,x)\,\, dx dt \\
=&\int_0^\infty  \int_\R \1_t(t) \Pi(u_1,u_2)(t,x)u_3(t,x)\,\, dx dt \\
=&\int_0^\infty  \int_\R  \Pi(\1_t(t)u_1,u_2)(t,x)u_3(t,x)\,\, dx dt \\
=&\int_0^\infty  \int_\R  \Pi(\1_{t,R}^{\text{high}}u_1,u_2)(t,x)u_3(t,x)\,\, dx dt  \\
&+ \int_0^\infty  \int_\R  \Pi(\1_{t,R}^{\text{low}}u_1,u_2)(t,x)u_3(t,x)\,\, dx dt \\
:=& \quad I_t^{\text{high}} \quad + \quad  I_t^{\text{low}}.
\end{align*} 
By Lemma \ref{Lemma 2.4} and Hölder's inequality, we can estimate $I_t^{\text{high}}$ like
\begin{align*}
|I_t^{\text{high}}|&= \left|\int_0^\infty  \int_\R  \Pi(\1_{t,R}^{\text{high}}u_1,u_2)(t,x)u_3(t,x)\,\, dx dt\right| \\
& \leq  \int_0^\infty \left| \int_\R  \Pi(\1_{t,R}^{\text{high}}u_1,u_2)(t,x)u_3(t,x)\,\, dx\right| dt \\
& \lesssim N_1^{1/2}   \int_0^\infty |\1_{t,R}^{\text{high}}| \| u_1(t)\|_{L^2_x} \| u_2(t)\|_{L^2_x}\| u_3(t)\|_{L^2_x}\,\,dt \\
&\leq N_1^{1/2}  \left(\int_0^\infty |\1_{t,R}^{\text{high}}| \| u_1(t)\|_{L^2_x}  \, dt\right) \, \| u_2\|_{L^\infty_tL^2_x}\| u_3\|_{L^\infty_tL^2_x} \\
& \leq  N_1^{1/2} \|\1_{t,R}^{\text{high}}\|_{L^q_t}\|u_1\|_{L^p_tL^2_x}  \| u_2\|_{L^\infty_tL^2_x}\| u_3\|_{L^\infty_tL^2_x} \\
& \leq N_1^{1/2} R^{\frac{1}{p}-1}\|u_1\|_{L^p_tL^2_x}  \| u_2\|_{L^\infty_tL^2_x}\| u_3\|_{L^\infty_tL^2_x}.
\end{align*}
Now, we turn our attention to the term $I_t^{\text{low}}$. Let $L \in 2^\Z$  be such that the condition
$$
L \gtrsim N_1N_2^2 +N_1N_2^{-2}
$$
does not hold. For $R \ll L,$ by Lemma \ref{lemma 2.2}, we have the decomposition
\begin{align*}
I_t^{\text{low}} =&\,\, I_\infty(Q_{\gtrsim L}(\1_{t,R}^{\text{low}}u_1),u_2,u_3) + I_\infty(Q_{\ll L}(\1_{t,R}^{\text{low}}u_1),Q_{\gg L}u_2,u_3) \\
& +I_\infty(Q_{\ll L}(\1_{t,R}^{\text{low}}u_1),Q_{\ll L}u_2,Q_{\gtrsim L}u_3)\\
 :=&\quad  I_1 + I_2 + I_3.
\end{align*}
By  Lemma \ref{Lemma 2.4} and Plancherel's identity we can estimate the contribution of $I_1$ and $I_2$  as
\begin{align*}
|I_1| &= |I_\infty(Q_{\gtrsim L}(\1_{t,R}^{\text{low}}u_1),u_2,u_3)| \\
%&=\left|\int_0^\infty  \int_\R  \Pi(Q_{\gtrsim L}(\1_{t,R}^{\text{low}}u_1),u_2)(t,x)u_3(t,x)\,\, dx \, dt\right| \\
%& \leq  \int_0^\infty \left| \int_\R  \Pi(Q_{\gtrsim L}(\1_{t,R}^{\text{low}}u_1),u_2)(t,x)u_3(t,x)\,\, dx\right| dt \\
&  \lesssim N_1^{1/2}   \int_0^\infty \| Q_{\gtrsim L}(\1_{t,R}^{\text{low}}u_1(t)\|_{L^2_x} \| u_2(t)\|_{L^2_x}\| u_3(t)\|_{L^2_x}\,\,dt \\
&  \lesssim N_1^{1/2}   \int_0^\infty \| Q_{\sim L}u_1(t)\|_{L^2_x} \| u_2(t)\|_{L^2_x}\| u_3(t)\|_{L^2_x}\,\,dt \\
&  \lesssim N_1^{1/2}  \left( \int_0^\infty \| Q_{\sim L}u_1(t)\|_{L^2_x} \| u_2(t)\|_{L^2_x}\,\,dt \right)\| u_3\|_{L_t^\infty L^2_x} \\
%& \lesssim  N_1^{1/2}  \| Q_{\sim L}u_1\|_{L^2_tL^2_x} \| u_2\|_{L^2_tL^2_x}\| u_3\|_{L_t^\infty L^2_x} \\
%& \lesssim  {N_1^{1/2}}L^{-1}\| \sigma(\tau, \xi) \tilde{u}_1\|_{L^2_tL^2_x} \| u_2\|_{L^2_tL^2_x}\| u_3\|_{L_t^\infty L^2_x} \\
&  \lesssim  {N_1^{1/2}}L^{-1} \|{u}_1\|_{X^{0,1}} \| u_2\|_{L^2_tL^2_x}\| u_3\|_{L_t^\infty L^2_x},
\end{align*}
and
\begin{align*}
	|I_2| &= |I_\infty(Q_{\ll L}(\1_{t,R}^{\text{low}}u_1),Q_\gtrsim u_2,u_3)| \\
	&=\left|\int_0^\infty  \int_\R  \Pi(Q_{\ll L}(\1_{t,R}^{\text{low}}u_1),Q_\gtrsim u_2)(t,x)u_3(t,x)\,\, dx \, dt\right| \\
	& \leq  \int_0^\infty \left| \int_\R  \Pi(Q_{\ll L}(\1_{t,R}^{\text{low}}u_1),Q_\gtrsim u_2)(t,x)u_3(t,x)\,\, dx\right| dt \\
	&  \lesssim N_1^{1/2}   \int_0^\infty \| Q_{\ll L}(\1_{t,R}^{\text{low}}u_1)(t)\|_{L^2_x} \| Q_\gtrsim u_2(t)\|_{L^2_x}\| u_3(t)\|_{L^2_x}\,\,dt \\
	&  \lesssim N_1^{1/2}   \int_0^\infty \| Q_{\sim L}u_1(t)\|_{L^2_x} \|Q_\gtrsim  u_2(t)\|_{L^2_x}\| u_3(t)\|_{L^2_x}\,\,dt \\
	%&  \lesssim N_1^{1/2} L^{-1}  \|u_1\|_{L^\infty_t L^2_x} \|Q_\gtrsim  u_2\|_{L_t^2 L^2_x}\| u_3\|_{L_t^2 L^2_x} \\
	& \lesssim N_1^{1/2} L^{-1}  \|u_1\|_{L^\infty_t L^2_x} \|u_2\|_{X^{0,1}}\| u_3\|_{L_t^2 L^2_x}.
\end{align*}
Finally, similar to what was done in $I_2$, we get
\begin{align*}
	|I_3|\lesssim N_1^{1/2} L^{-1}  \|u_1\|_{L^\infty_t L^2_x} \| u_2\|_{L_t^2 L^2_x}  \|u_3\|_{X^{0,1}}.
\end{align*}
\qed
\begin{remark}
 The case of $\bar{Q}$ follows in an entirely analogous manner.
\end{remark}

For $s \in \R$ we define our resolution space being
$$
\M^s = \M^s(\R \times \R):= L^\infty_t Z^s_x \cap X^{s-1,1}
$$
endowed with the norm
$$
\|\cdot \|_{\M^s} : = \|\cdot \|_{L^\infty_t Z^s_x} +\|\cdot \|_{X^{s-1,1}} .
$$
Fix $T<1$. For working satisfactorily in our resolution space we have to extend a function $f: (0,T) \times \R \rightarrow \R$ from $(0,T)$ to all real line $\R$. Following the work of Masmoudi and Nakanishi \cite{MR2659158} we introduce the extension operator $\rho_T$ defined by
$$
\rho_T [f](t,x):= \eta(t)f\left(T \mu\left(t \cdot T^{-1}\right), x \right)
$$
where $\mu(t) = \max(1-|t-1|,0)$ and cut-off the function $\eta$ was introduced in Section \ref{seção das notação}. The operator $\rho_T[f]$ is defined in $\R \times \R$ and  $\rho_T[f](t,\cdot) \equiv f(\cdot)$ whenever $t \in (0,T)$. Furthermore we have the following properties:
\begin{itemize}
\item [-] The operator $\rho_T: X^{s,b}_{T} \rightarrow  X^{s,b}$ is bounded for all $b \leq 1$ and $s \in \R$;
\item [-] The operator $\rho_T: L^p_T(B) \rightarrow  L^p(B) $ is bounded for all $1 \leq p \leq \infty$ and $B$ a Banach space,
\end{itemize}
and these two properties are uniformly for $0<T<1$ (See \cite{MR2659158} for more details). Using this approach, we establish the first \textit{a priori} estimate, ensuring that solutions of \eqref{BO} in $ L^\infty_T Z^s_x $ also lie in $ \M_T^s $. Moreover, the difference between two solutions in $ L^\infty_T Z^{s-1}_x $ belongs to $ \M_T^{s-1} $.

\begin{lemma}\label{lemma 3.1}
Let $T \in (0,1)$, $s> 1/2$ and $u \in L^\infty_T Z^s_x$ be a solution to the equation \eqref{BO2}  with initial data $u_0 \in Z^s_x(\R)$. Then $u \in \M^s_T$ and
\begin{equation}
\|u\|_{\M^s_T} \lesssim \left(1+\|u\|_{L^\infty_T Z^{\frac{1}{2}+}_x}\right) \|u\|_{L^\infty_T Z^s_x}.
\end{equation}
Moreover, for any pair $u,v \in L^\infty_TZ^s_x$ of solutions to IVP \eqref{BO2}  associated with initial data $u_0, v_0 \in Z^s_x(\R)$ respectively, the following inequalities holds
\begin{equation}
	\| u -v \|_{\M^{s}_T} \lesssim \|u-v\|_{L^\infty_T Z^{s}_x} + \|u+v\|_{L^\infty_T Z^s_x}\|u-v\|_{L^\infty_T Z^{s}_x},
\end{equation}
and
\begin{equation}
\| u -v \|_{\M^{s-1}_T} \lesssim \|u-v\|_{L^\infty_T Z^{s-1}_x} + \|u+v\|_{L^\infty_T Z^s_x}\|u-v\|_{L^\infty_T Z^{s-1}_x}.
\end{equation}
\end{lemma}
\dem  Using the extension operator, we define the extension $\textbf{u}(t) := \rho_T[u](t)$, which is defined for all $t \in \R$ but supported on the interval $(0, T)$. For simplicity, we will henceforth omit the boldface in the notation. 

By the properties of the extension operator, it suffices to estimate the $X^{s,b}_{T}$-norm of $u$. We start by expressing the solution through its Duhamel formulation
$$
u(t)=  V(t)u_0 +  \int_0^t V(t-s)u \partial_x u(s)\, ds.
$$
Then, by the energy estimates in Bourgain spaces, we obtain for $s> \frac{1}{2}$
\begin{align*}
\|u\|_{X_T^{s-1,1}} &\lesssim \|u_0\|_{Z^{s-1}_x} + \|u \partial_x u \|_{X_T^{s-1,0}} \\
%&= \|u_0\|_{Z^{s-1}_x} + \|\partial_x( u^2) \|_{X^{s-1,0}} \\
& \lesssim  \|u_0\|_{Z^s_x} + \|u^2 \|_{X_T^{s,0}} \\
%& \lesssim\|u_0\|_{Z^s_x} + \|u^2 \|_{L^2_t Z^s_x} \\
& \lesssim \|u_0\|_{Z^s_x} + \|u^2 \|_{L^\infty_t H^s_x} \\
%& \lesssim \|u_0\|_{Z^s_x} + \|u \|_{L^\infty_t Z^{\frac{1}{2}+}_x}\|u \|_{L^\infty_t Z^s_x}\\
&\lesssim\|u_0\|_{Z^s_x} + \|u \|_{L^\infty_t H^{\frac{1}{2}+}_x}\|u \|_{L^\infty_t H^s_x}\\
%&\lesssim\|u_0\|_{Z^s_x} + \|u \|_{L^\infty_t Z^{\frac{1}{2}+}_x}\|u \|_{L^\infty_t Z^s_x}\\
%&\lesssim \sup_{t \in [0,T]}\|u(t)\|_{Z^s_x} + \|u \|_{L^\infty_t Z^{\frac{1}{2}+}_x}\|u \|_{L^\infty_t Z^s_x}\\
& \lesssim \|u\|_{L^\infty_TZ^s_x} + \|u \|_{L^\infty_t Z^{\frac{1}{2}+}_x}\|u \|_{L^\infty_t Z^s_x}\\
%& \lesssim \|u\|_{L^\infty_TZ^s_x} + \|u \|_{L^\infty_t Z^{\frac{1}{2}+}_x}\|u \|_{L^\infty_t Z^s_x}\\
& \lesssim \left(1+\|u\|_{L^\infty_T Z^{\frac{1}{2}+}_x}\right) \|u\|_{L^\infty_T Z^s_x}.
\end{align*}
Similarly, if  $u,v \in L^\infty_TZ^s_x$ are solutions of \eqref{BO2}  with initial data $u_0, v_0 \in Z^s_x(\R)$ respectively we have
\begin{align*}
	\|u-v\|_{X_T^{s-2,1}} &\lesssim \|u_0-v_0\|_{Z^{s-2}_x} + \|u^2 - v^2\|_{L^2_tZ^{s-1}_x }\\
%&= \|u_0-v_0\|_{Z^{s-2}_x} + \|(u+v)(u-v)\|_{L^2_tZ^{s-1}_x }\\
& \lesssim \|u_0-v_0\|_{Z^{s-1}_x} + \|(u+v)(u-v)\|_{L^\infty_TZ^{s-1}_x }\\
&\lesssim \|u-v\|_{L^\infty_T Z^{s-1}_x} + \|u+v\|_{L^\infty_TZ^{s}_x } \|u-v\|_{L^\infty_TZ^{s-1}_x },
\end{align*}
where in the last inequality we use Lemma \ref{lemma do sobolev} with $s_1=s$ and $s_2=s_3=s-1$ for $s>1/2$. The case  $\|u-v\|_{X_T^{s-1,1}}$ is analogous. \qed

To proceed, we present the fundamental energy estimate that will serve as the foundation for the compactness argument employed to ensure the existence of solutions to \eqref{BO2}.

\begin{theorem}[Energy Estimate]\label{des de energia}
	Let $T \in (0,1)$, $s> 1/2$ and $u \in L^\infty_T Z^s_x$ be a solution to the IVP \eqref{BO2} with initial data $u_0 \in Z^s_x(\R)$. Then we have the following energy estimate
	\begin{equation}
		\|u\|^2_{L_T^\infty Z^s_x} \lesssim \|u_0\|^2_{Z^s_x} + \|u\|_{L^\infty_T L^2_x} \|u\|_{L^\infty_T Z^s_x}^2 + \left(1+\|u\|_{L^\infty_T Z^{\frac{1}{2}+}_x}\right)^2 \|u\|_{L^\infty_T Z^{\frac{1}{2}+}_x} \|u\|_{L^\infty_T Z^s_x}^2 + \|u\|_{L^\infty_T Z^s_x}^3.
	\end{equation}
\end{theorem}
\dem We begin by analyzing, formally, the classical estimate for the dyadic component $ P_N u $, which is localized around the spatial frequency $ N \in 2^\Z $. Applying the Littlewood-Paley operator $ P_N $ to the equation in \eqref{BO2}, we take the $ L^2_x $-inner product with $ P_N u $, multiply by $ \langle N \rangle^{2s} $ (for $ s > \frac{1}{2} $), and integrate over the time interval $(0, t)$, with $ t \in (0, T) $. This procedure leads us to the following identity:
\begin{equation}\label{expression1}
    \langle N\rangle^{2s} \|P_N u(t)\|_{L^2_x} = \langle N\rangle^{2s} \|P_N u(0)\|_{L^2_x} + \langle N\rangle^{2s}\int_0^t \int_\R \partial_x P_N(u^2)(t)P_N u (t) \, dx \, dt.
\end{equation}
Similarly, applying the operators $P_N$ and $\partial_x^{-1}$ in the equation \eqref{BO2}, taking the $L^2_x$-inner product with $\partial_x^{-1}u$ and integrating on $(0,t)$ with $ t \in (0,T)$ we obtain
\begin{equation}\label{expression2}
\|P_N \partial_x^{-1}u(t)\|_{L^2_x} =\|P_N \partial_x^{-1}u(0)\|_{L^2_x} + \int_0^t \int_\R P_N(u^2)(t)P_N \partial_x^{-1}u (t) \, dx \, dt.
\end{equation}
By \eqref{expression1}, \eqref{expression2}, Bernstein's inequalities and integrating by parts we have
\begin{equation}\label{3-11}
\begin{aligned}
    \| P_N u \|_{L^\infty_TZ^s_x}^2 \lesssim \| P_N u(0) \|_{Z^s_x}^2 &+ \sup_{t \in (0,T)} \langle N\rangle^{2s} \left|\int_0^t \int_\R P_N(u^2)(t) \partial_x P_N u (t) \, dx \, dt \right|\\
    &+ \sup_{t \in (0,T)}  \left|\int_0^t \int_\R P_N(u^2)(t)  P_N \partial_x^{-1}u (t) \, dx \, dt \right|. 
\end{aligned}
\end{equation}
We aim to sum over all $ N \in 2^\Z $ in inequality \eqref{3-11} to apply the Theorem \ref{Littlewood-Paley Theorem} and derive an estimate of the form \eqref{des de energia}. To achieve this, we first need to obtain suitable bounds for the following quantities: 
\begin{equation}
I: =  \sum_{N \in 2^\Z}\sup_{t \in (0,T)} \langle N\rangle^{2s}	 \left|\int_0^t \int_\R P_N(u^2)(t) \partial_x P_N u (t) \, dx \, dt \right|.
\end{equation}
and
\begin{equation}\label{rotatingpart}
    J: =  \sum_{N \in 2^\Z}\sup_{t \in (0,T)} \left|\int_0^t \int_\R P_N(u^2)(t) P_N \partial_x^{-1}u (t) \, dx \, dt \right|.
\end{equation}
Denoting $ u_N(t) := P_N u(t) $ (and other possible variations, such as $ u_{\gg N}, u_{\sim N}, \text{etc} $). From this point onward, we also consider an extension $\mathbf{u}$ of $u$ with temporal support in $(-1,1)$ such that
\[
\|\mathbf{u}\|_{\mathcal{M}^s} \lesssim \|u\|_{\mathcal{M}_T^s}.
\]  
Since there is no risk of confusion, we will simply denote $\mathbf{u}$ by $u$. We begin by estimating $ I $. First, observe that
\begin{equation}\label{expressão do Pn 1}
P_N(u^2) = P_N(u_{\gtrsim N} u_{\gtrsim N}) +2 P_N(u_{\ll N} u).
\end{equation}
Using the Taylor expansion with the integral remainder, we have for $ \xi_1 \in \mathbb{R} $
$$
\varphi_N(\xi)=\varphi_N(\xi - \xi_1) - N^{-1}\xi_1 \int_0^1 (1- h) \varphi_N(\xi - h \xi_1)\, dh,
$$
where $\varphi_N$ is the Littlewood-Paley projector defined in Section \ref{seção das notação}. Therefore
\begin{align*}
 \widehat{P_N(u_{\ll N} u)} (\xi) &= \varphi_N (\xi) (\widehat{u_{\ll N}} * \widehat{u}) (\xi) \\
 &= \int_{\mathbb{R}} \varphi_N (\xi) \widehat{u_{\ll N}} (\xi_1) \widehat{u}(\xi-\xi_1) \, d\xi_1\\
 & = \int_{\mathbb{R}}  \widehat{u_{\ll N}} (\xi_1) \varphi_N(\xi - \xi_1)\widehat{u}(\xi-\xi_1) \, d\xi_1\\
 &\quad + \int_{\mathbb{R}}  \xi_1 \widehat{u_{\ll N}} (\xi_1) \widehat{u}(\xi-\xi_1) \cdot \left[- N^{-1} \int_0^1 (1- h) \varphi_N^\prime(\xi - h \xi_1)\, dh \right] \, d\xi_1\\
& = \int_{\mathbb{R}}  \widehat{u_{\ll N}} (\xi_1) \widehat{P_N u}(\xi-\xi_1) \, d\xi_1\\
 &\quad + \int_{\mathbb{R}}   \widehat{\partial_x u_{\ll N}} (\xi_1) \widehat{u}(\xi-\xi_1) \cdot \left[-i N^{-1} \int_0^1 (1- h) \varphi_N^\prime (\xi - h \xi_1)\, dh \right] \, d\xi_1\\
 & = \widehat{u_{\ll N} P_Nu} (\xi) +\int_{\mathbb{R}}   \widehat{\partial_x u_{\ll N}} (\xi_1) \widehat{u}(\xi-\xi_1) \cdot \Psi(\xi,\xi_1) \, d\xi_1.
\end{align*}
Thus, we obtain
\begin{equation}\label{expressao do Pn 2}
P_N(u_{\ll N} u) = u_{\ll N} P_Nu + \Pi_\Psi (\partial_x u_{\ll N}, u),
\end{equation}
with
$$
 \Psi(\xi,\xi_1):=-i N^{-1} \int_0^1 (1- h) \varphi_N^\prime(\xi - h \xi_1)\, dh \in L^\infty(\R^2).
$$
By substituting \eqref{expressao do Pn 2} into \eqref{expressão do Pn 1}, we obtain that
\begin{equation}\label{expressao do Pn 3}
P_N(u^2) = P_N(u_{\gtrsim N} u_{\gtrsim N}) +2  u_{\ll N} P_Nu + 2 \Pi (\partial_x u_{\ll N}, u)
\end{equation}
 Substituting \eqref{expressao do Pn 3} into $ I $, integrating by parts, and recalling that the definition of the operators $ I_t $ is given by  
\[
I_t(u_1, u_2, u_3) = \int_0^t \int_\mathbb{R} \Pi(u_1, u_2)(t, x) u_3(t, x) \, dx \, dt,
\]  
we obtain (see Proposition 3.4 in \cite{MR3397003} for more details)
\begin{align*}
I &\lesssim \sum_{N \in 2^\Z}\sum_{N_1 \gtrsim N} N \langle N_1 \rangle^{2s} \sup_{t \in (0,T)} |I_t(u_N,u_{\sim N_1}, u_{N_1})|
\end{align*}
for some pseudoproduct $\Pi=\Pi_\Psi$. Rewrite $I$ as being $I=I_{N \lesssim 1} + I_{N \gg 1}$ with
$$
I_{N \lesssim 1}:=\sum_{N \lesssim 1}\sum_{N_1 \gtrsim N} N \langle N_1 \rangle^{2s} \sup_{t \in (0,T)} |I_t(u_N,u_{\sim N_1}, u_{N_1})|
$$
and
$$
I_{N \gg 1}:=\sum_{N \gg 1}\sum_{N_1 \gtrsim N} N \langle N_1 \rangle^{2s} \sup_{t \in (0,T)} |I_t(u_N,u_{\sim N_1}, u_{N_1})|.
$$
We will start by estimating the term $I_{N \lesssim 1}$. Then, by the first part of Lemma \ref{lemma 2.2}
\begin{align*}
I_{N \lesssim 1}&\lesssim \sum_{N \lesssim 1}\sum_{N_1 \gtrsim N} N \langle N_1 \rangle^{2s} N^{1/2}\sup_{t \in (0,T)}\|u_N\|_{L^\infty_t L^2_x}\|u_{\sim N_1}\|_{L^2_tL^2_x}\|u_{N_1}\|_{L^2_tL^2_x}\\
&\lesssim \sum_{N \lesssim 1}\sum_{N_1 \gtrsim N} \langle N_1 \rangle^{2s} \|u_N\|_{L^\infty_T L^2_x}\|u_{ N_1}\|_{L^2_TL^2_x}\|u_{N_1}\|_{L^2_TL^2_x}\\
&=\sum_{N \lesssim 1}\sum_{N_1 \gtrsim N} \langle N_1 \rangle^{2s} \|u_N\|_{L^\infty_T L^2_x}\|u_{ N_1}\|_{L^2_TL^2_x}^2\\
&=\sum_{N \lesssim 1}\sum_{N_1 \gtrsim N} \|u_N\|_{L^\infty_T L^2_x}\|u_{ N_1}\|_{L^\infty_T H^s_x}^2\\
&=\sum_{N \lesssim 1}\sum_{N_1 \gtrsim N} \|u_N\|_{L^\infty_T L^2_x}\|u_{ N_1}\|_{L^\infty_T Z^s_x}^2\\
& \lesssim  \|u\|_{L^\infty_T L^2_x} \|u\|_{L^\infty_T Z^s_x}^2.
\end{align*}
Now we estimate the term $I_{N \gg 1}$. By the second inequality in Lemma \ref{l3.4} 
\begin{align*}
I_{N \gg 1} &\lesssim \sum_{N \gg 1}\sum_{N_1 \gtrsim N} N \langle N_1 \rangle^{2s}\Big[
N^{\frac{1}{2}}R^{\frac{1}{p}-1}\|u_N\|_{L^p_T L^2_x}\|u_{ N_1}\|_{L^\infty_TL^2_x}^2\\
&\qquad \qquad +N^{\frac{1}{2}}L^{-1} \|u_N\|_{X_T^{0,1}}\|u_{ N_1}\|_{L^2_TL^2_x}\|u_{ N_1}\|_{L^\infty_TL^2_x} \\
&\qquad \qquad +N^{\frac{1}{2}}L^{-1}\|u_{ N}\|_{L^\infty_TL^2_x} \|u_{N_1}\|_{X_T^{0,1}}  \|u_{ N_1}\|_{L^2_TL^2_x}\Big]\\
&\lesssim \sum_{N \gg 1}\sum_{N_1 \gtrsim N} N^{\frac{3}{2}}\langle N_1 \rangle^{2s}R^{\frac{1}{p}-1}\|u_N\|_{L^p_T L^2_x}\|u_{ N_1}\|_{L^\infty_TL^2_x}^2 \\
&\quad +\sum_{N \gg 1}\sum_{N_1 \gtrsim N} N^{\frac{3}{2}}\langle N_1 \rangle^{2s}L^{-1} \|u_N\|_{X_T^{0,1}}\|u_{ N_1}\|_{L^2_TL^2_x}\|u_{ N_1}\|_{L^\infty_TL^2_x}  \\
&\quad +\sum_{N \gg 1}\sum_{N_1 \gtrsim N} N^{\frac{3}{2}}\langle N_1 \rangle^{2s}L^{-1} \|u_{ N}\|_{L^\infty_TL^2_x} \|u_{N_1}\|_{X_T^{0,1}}  \|u_{ N_1}\|_{L^2_TL^2_x}\\
&\lesssim \sum_{N \gg 1}\sum_{N_1 \gtrsim N} N^{\frac{3}{2}}\langle N_1 \rangle^{2s}R^{\frac{1}{p}-1}\|u_N\|_{L^p_T L^2_x}\|u_{ N_1}\|_{L^\infty_TL^2_x}^2 \\
&\quad +\sum_{N \gg 1}\sum_{N_1 \gtrsim N} N^{\frac{3}{2}}\langle N_1 \rangle^{2s}L^{-1} \|u_N\|_{X_T^{0,1}}\|u_{ N_1}\|_{L^\infty_TL^2_x}^2\\
&\quad +\sum_{N \gg 1}\sum_{N_1 \gtrsim N} N^{\frac{3}{2}}\langle N_1 \rangle^{2s}L^{-1} \|u_{ N}\|_{L^\infty_TL^2_x} \|u_{N_1}\|_{X_T^{0,1}}  \|u_{ N_1}\|_{L^2_TL^2_x}\\
&\lesssim \sum_{N \gg 1}\sum_{N_1 \gtrsim N} N^{\frac{3}{2}}R^{\frac{1}{p}-1}\|u_N\|_{L^p_T L^2_x}\|u_{ N_1}\|_{L^\infty_TH^s_x}^2 \\
&+\sum_{N \gg 1}\sum_{N_1 \gtrsim N} N^{\frac{3}{2}}L^{-1} \|u_N\|_{X_T^{0,1}}\|u_{ N_1}\|_{L^\infty_TH^s_x}^2  \\
&\quad +\sum_{N \gg 1}\sum_{N_1 \gtrsim N} N^{\frac{3}{2}}L^{-1} \|u_{ N}\|_{L^\infty_TL^2_x} \|u_{N_1}\|_{X_T^{s,1}}  \|u_{ N_1}\|_{L^2_TH^s_x} \\
\end{align*}
\begin{align*}
&\lesssim \sum_{N \gg 1}\sum_{N_1 \gtrsim N} N^{\frac{3}{2}}R^{\frac{1}{p}-1}\|u_N\|_{L^p_T L^2_x}\|u_{ N_1}\|_{L^\infty_TZ^s_x}^2 \\
&\quad +\sum_{N \gg 1}\sum_{N_1 \gtrsim N} N^{\frac{3}{2}}L^{-1} \|u_N\|_{X_T^{0,1}}\|u_{ N_1}\|_{L^\infty_TZ^s_x}^2  \\
&\quad +\sum_{N \gg 1}\sum_{N_1 \gtrsim N} N^{\frac{3}{2}}L^{-1} \|u_{ N}\|_{L^\infty_TL^2_x} \|u_{N_1}\|_{X_T^{s,1}}  \|u_{ N_1}\|_{L^2_TZ^s_x}.
\end{align*}
Choosing $p= +\infty$, $R = N^{3/2}N_1^{1/8}$, $L=NN_1^2$ and remembering that   $\langle N \rangle^{2s} \gg 1$  we can obtain
\begin{align*}
I_{N \gg 1}  \lesssim \|u\|_{\M_T^{\frac{1}{2}+}}\|u\|_{\M_T^{s}}\|u\|_{L^2_TZ^s_x}.
\end{align*}
Gathering the estimates for $I_{N \lesssim 1}$ and $I_{N \gg 1}$ and using Lemma \ref{lemma 3.1}  we obtain
\begin{align}\label{estimativa do I}
\nonumber I &\lesssim \|u\|_{L^\infty_T L^2_x} \|u\|_{L^\infty_T Z^s_x}^2 + \|u\|_{\M_T^{\frac{1}{2}+}}\|u\|_{\M_T^{s}}\|u\|_{L^2_TZ^s_x}\\
& \lesssim  \|u\|_{L^\infty_T L^2_x} \|u\|_{L^\infty_T Z^s_x}^2 + \left(1+\|u\|_{L^\infty_T Z^{\frac{1}{2}+}_x}\right)^2 \|u\|_{L^\infty_T Z^{\frac{1}{2}+}_x} \|u\|_{L^\infty_T Z^s_x}^2.
\end{align}
Now, to deal with the term $J$, in \eqref{rotatingpart}, observe that, by the Cauchy-Schwarz inequality we have
\begin{align*}
 \left| \int_0^t\int_\R P_N(u^2)(t) P_N \partial_x^{-1}u (t) \, dx \, dt\right|& \leq \int_0^t \int_\R |P_N(u^2)(t)|| P_N \partial_x^{-1}u (t)| \, dx \, dt\\
 &\leq \int_0^t\| P_N(u^2)(t)\|_{L^2_x}\|P_N \partial_x^{-1}u (t)\|_{L^2_x}\, dt\\
 &\lesssim \| P_N(u^2)\|_{L^\infty_tL^2_x}\|P_N \partial_x^{-1}u\|_{L^\infty_tL^2_x}.
\end{align*}
Therefore, by the Sobolev embedding and Cauchy-Schwarz in $l^2(\mathbb{Z})$ we have, for $s>1/2$
\begin{align*}
J &= \sum_{N \in 2^\Z}\sup_{t \in (0,T)}	 \left|\int_0^t \int_\R P_N(u^2)(t) P_N \partial_x^{-1}u (t) \, dx \, dt \right| \\
&\lesssim \sum_{N \in 2^\Z}  \| P_N(u^2)\|_{L^\infty_TL^2_x}\|P_N \partial_x^{-1}u\|_{L^\infty_TL^2_x}\\
%&= \sum_{N \in 2^\Z}  \left(\| P_N(u^2)\|_{L^\infty_TL^2_x} \right)\left(\|P_N \partial_x^{-1}u\|_{L^\infty_TL^2_x}\right)\\
&\lesssim \left(  \sum_{N \in 2^\Z}  \| P_N(u^2)\|_{L^\infty_TL^2_x}^2\right)^{1/2}\left(  \sum_{N \in 2^\Z}  \|P_N \partial_x^{-1}u\|_{L^\infty_TL^2_x}^2\right)^{1/2}\\
&\lesssim \left(  \sum_{N \in 2^\Z}  \langle N\rangle^{2s}\| P_N(u^2)\|_{L^\infty_TL^2_x}^2\right)^{1/2}\left(  \sum_{N \in 2^\Z}  \|P_N \partial_x^{-1}u\|_{L^\infty_TL^2_x}^2\right)^{1/2}\\
& \lesssim \|u^2\|_{L^\infty_TH^s_x} \|\partial_x^{-1}u\|_{L^\infty_T L^2_x} \\
&\lesssim \|u\|_{L^\infty_TH^s_x}^2\|\partial_x^{-1}u\|_{L^\infty_T H^s_x}\\
&\lesssim \|u\|_{L^\infty_TZ^s_x}^3.
\end{align*}
Finally, summing over $N\in 2^\Z$ in inequality \eqref{3-11} and applying the Theorem \ref{Littlewood-Paley Theorem} lead to
$$
\|u\|^2_{L_T^\infty Z^s_x} \lesssim \|u_0\|^2_{Z^s_x} + \|u\|_{L^\infty_T L^2_x} \|u\|_{L^\infty_T Z^s_x}^2 + \left(1+\|u\|_{L^\infty_T Z^{\frac{1}{2}+}_x}\right)^2 \|u\|_{L^\infty_T Z^{\frac{1}{2}+}_x} \|u\|_{L^\infty_T Z^s_x}^2 + \|u\|_{L^\infty_T Z^s_x}^3.
$$
\qed

It is also crucial to estimate the difference between two solutions, as indicated by the following theorem.
\begin{theorem}\label{diferença em Zs}
Let $T \in (0,1)$ and $u_1,u_2 \in L^\infty_TZ^s_x$ be two solutions of \eqref{BO2} associated to $u_1^0, u_2^0 \in Z^s_x(\R)$ respectively, where $s>\frac{1}{2}$. Then
\begin{equation}
\|u_1-u_2\|^2_{L^\infty_T Z^{s-1}_x} \lesssim \|u_1^0-u_2^0\|^2_{Z^{s-1}_x} + \|u_1 + u_2\|_{\M_T^{s}}\|u_1 - u_2\|^2_{\M_T^{s-1}}.
\end{equation}
\end{theorem}
\dem We start by defining the differences $ u_0 := u_1^0 - u_2^0 $, $ u := u_1 - u_2 $, and $ v := u_1 + u_2 $. Observe that $ u $ satisfies the following IVP
\begin{equation}\label{equação da diferença}
	\left\{
	\begin{array}{l}
		\partial_t u + \mathcal{H} \partial_x^2 u + \partial_x(uv)=  \gamma\partial_x^{-1} u, \quad (t,x)  \in \R \times \R\\
		u(0,x)=u_0(x) \quad x \in \R.
	\end{array}
	\right.
\end{equation}

By following the analogous procedure outlined in the proof of Theorem \ref{des de energia}, we can derive
\begin{align}\label{desigualdade u-v}
\nonumber\| P_N u \|_{L^\infty_TZ^{s-1}_x}^2 \lesssim \| P_N u_0 \|_{Z^{s-1}_x}^2 &+ \sup_{t \in (0,T)} \langle N\rangle^{2(s-1)} \left|\int_0^t \int_\R P_N(uv)(t) \partial_x P_N u (t) \, dx \, dt \right|\\
&+ \sup_{t \in (0,T)}\left|\int_0^t \int_\R P_N(uv)(t) \partial_x^{-1} P_N u (t) \, dx \, dt \right|.
\end{align}
Using a similar approach to the previous result, our goal is to sum over all dyadic numbers $ N \in 2^\mathbb{Z} $ and apply the Theorem \ref{Littlewood-Paley Theorem} to inequality \eqref{desigualdade u-v}. To accomplish this, we first require an estimate for:
$$
I := \sum_{N \in 2^\Z}  \sup_{t \in (0,T)} \langle N\rangle^{2(s-1)} \left|\int_0^t \int_\R P_N(uv)(t) \partial_x P_N u (t) \, dx \, dt \right|,
$$
and
$$
J := \sum_{N \in 2^\Z}  \sup_{t \in (0,T)} \left|\int_0^t \int_\R P_N(uv)(t) \partial_x^{-1}P_N u (t) \, dx \, dt \right|.
$$
Following the same steps performed in the proof of Theorem \ref{des de energia}  we can decompose $I$ as being (see Proposition 3.5 in \cite{MR3397003} for more details)
\begin{align*}
I &\lesssim \sum_{N \in 2^\Z}\sum_{N_1 \gtrsim N} N \langle N_1 \rangle^{2(s-1)} \sup_{t \in (0,T)} |I_t(v_N,u_{\sim N_1}, u_{N_1})| \\
&\quad +\sum_{N \in 2^\Z}\sum_{N_1 \gtrsim N} N_1 \langle N_1 \rangle^{2(s-1)} \sup_{t \in (0,T)} |I_t(v_{\sim N_1},u_{N}, u_{N_1})| \\
&\quad +\sum_{N \in 2^\Z}\sum_{N_1 \gtrsim N} N \langle N \rangle^{2(s-1)} \sup_{t \in (0,T)} |I_t(v_{N_1},u_{ N_1}, u_{N})| \\
&:= I_1 + I_2 + I_3.
\end{align*}
We will analyze the terms $I_j$ separately. Note that, since we have $s> \frac{1}{2}$
\begin{align*}
I_3 &=\sum_{N \in 2^\Z}\sum_{N_1 \gtrsim N} N \langle N \rangle^{2(s-1)} \sup_{t \in (0,T)} |I_t(v_{N_1},u_{ N_1}, u_{N})|\\
&\lesssim \sum_{N \in 2^\Z}\sum_{N_1 \gtrsim N} N_1 \langle N_1 \rangle^{2(s-1)} \sup_{t \in (0,T)} |I_t(v_{\sim N_1},u_{N}, u_{N_1})| \\
&= I_2.
\end{align*}
Therefore, it suffices to estimate $ I_1 $ and $ I_2 $. As in the proof of Theorem \ref{des de energia}, we rewrite $ I_j $ as
$$
I_j = I^{(j)}_{N \lesssim 1} + I^{(j)}_{N \gg 1} \quad j=1,2,
$$
as being the contribution of the sum over $N \lesssim 1$ and $N \gg 1$ respectively. By the first inequality in Lemma \ref{lemma 2.2}
\begin{align*}
I^{(1)}_{N \lesssim 1} &= \sum_{N \lesssim 1}\sum_{N_1 \gtrsim N} N \langle N_1 \rangle^{2(s-1)} \sup_{t \in (0,T)} |I_t(v_N,u_{\sim N_1}, u_{N_1})|\\
& \lesssim \sum_{N \lesssim 1}\sum_{N_1 \gtrsim N} N \langle N_1 \rangle^{2(s-1)} N^{1/2} \|v_N\|_{L^\infty_T L^2_x}\|u_{N_1}\|_{L^2_TL^2_x}^2\\
& \lesssim \sum_{N \lesssim 1}\sum_{N_1 \gtrsim N}  N^{3/2} \|v_N\|_{L^\infty_T L^2_x}\|u_{N_1}\|_{L^2_TH^{s-1}_x}^2 \\
&\lesssim \sum_{N \lesssim 1}\sum_{N_1 \gtrsim N} \|v_N\|_{L^\infty_T L^2_x}\|u_{N_1}\|_{L^2_TH^{s-1}_x}^2 \\
& \lesssim \|v\|_{L^\infty_T L^2_x}\|u\|_{L^2_TH^{s-1}_x}^2\\
& \lesssim \|v\|_{L^\infty_T L^2_x}\|u\|_{L^2_TZ^{s-1}_x}^2,
\end{align*}
since $N \lesssim 1$. For the $I^{(2)}_{N \lesssim 1}$
\begin{align*}
I^{(2)}_{N \lesssim 1} &\lesssim \sum_{N \lesssim 1}\sum_{N_1 \gtrsim N} N_1 \langle N_1 \rangle^{2(s-1)} \sup_{t \in (0,T)} |I_t(v_{\sim N_1},u_{N}, u_{N_1})| \\
& \lesssim \sum_{N \lesssim 1}\sum_{N_1 \gtrsim N} N_1 \langle N_1 \rangle^{2(s-1)} N^{1/2} \|v_{N_1}\|_{L^2_T L^2_x}\|u_{N}\|_{L^2_TL^2_x}\|u_{N_1}\|_{L^\infty_TL^2_x}\\
&\lesssim \sum_{N \lesssim 1}\sum_{N_1 \gtrsim N} N_1 \langle N_1 \rangle^{s-1} \|v_{N_1}\|_{L^2_T L^2_x}\|u_{N}\|_{L^2_TL^2_x}\|u_{N_1}\|_{L^\infty_TH^{s-1}_x}\\
&\lesssim \sum_{N \lesssim 1}\sum_{N_1 \gtrsim N} N_1 \langle N_1 \rangle^{-1}  \|v_{N_1}\|_{L^2_T H^s_x}\|u_{N}\|_{L^2_TL^2_x}\|u_{N_1}\|_{L^\infty_TH^{s-1}_x}\\
&\lesssim \sum_{N \lesssim 1}\sum_{N_1 \gtrsim N}  \|v_{N_1}\|_{L^2_T H^s_x}\|u_{N}\|_{L^2_TH^{-\frac{1}{2}}_x}\|u_{N_1}\|_{L^\infty_TH^{s-1}_x}\\
&\lesssim \sum_{N \lesssim 1}\sum_{N_1 \gtrsim N}  \|v_{N_1}\|_{L^\infty_T H^s_x}\|u_{N}\|_{L^\infty_TH^{-\frac{1}{2}}_x}\|u_{N_1}\|_{L^\infty_TH^{s-1}_x}\\
& \lesssim  \|v\|_{L^\infty_T H^s_x}\|u\|_{L^\infty_TH^{-\frac{1}{2}}_x}\|u\|_{L^\infty_TH^{s-1}_x}\\
& \lesssim  \|v\|_{L^\infty_T Z^s_x}\|u\|_{L^\infty_TZ^{-\frac{1}{2}}_x}\|u\|_{L^\infty_TZ^{s-1}_x}.
\end{align*}
Now, we estimate the terms $ I^{(1)}_{N \gg 1} $ and $ I^{(2)}_{N \gg 1} $. Using the second inequality of Lemma \ref{lemma 2.2} with $ p = +\infty $, $ R = N^{3/2}N_1^{1/8} $, and $ L = NN_1^2 $, we obtain

\begin{align*}
I^{(1)}_{N \gg 1} &= \sum_{N \gg 1}\sum_{N_1 \gtrsim N} N \langle N_1 \rangle^{2(s-1)} \sup_{t \in (0,T)} |I_t(v_N,u_{\sim N_1}, u_{N_1})|  \\
&\lesssim  \sum_{N \gg 1}\sum_{N_1 \gtrsim N} N \langle N_1 \rangle^{2(s-1)} \Big(  N^{-1} N_1^{-\frac{1}{8}}\|v_N\|_{L^\infty_TL^2_x}  \| u_{N_1}\|_{L^\infty_TL^2_x}^2\\
& \quad  + {N^{-1/2}}N_1^{-2} \|{v}_N\|_{X_T^{0,1}} \| u_{N_1}\|_{L^2_TL^2_x}\| u_{N_1}\|_{L_T^\infty L^2_x} \\
&\quad + N^{-1/2} N_1^{-2}  \|v_N\|_{L^\infty_T L^2_x}\|u_{N_1}\|_{X_T^{0,1}}\| u_{N_1}\|_{L_T^2 L^2_x} \Big)\\
&\lesssim  \sum_{N \gg 1}\sum_{N_1 \gtrsim N} N \langle N_1 \rangle^{2(s-1)} \Big(  N^{-1} N_1^{-\frac{1}{8}}\|v_N\|_{L^\infty_TL^2_x}  \| u_{N_1}\|_{L^\infty_TL^2_x}^2\\
& \quad  + {N^{-1/2}}N_1^{-2} \|{v}_N\|_{X_T^{0,1}}\| u_{N_1}\|_{L_T^\infty L^2_x}^2 \\
&\quad + N^{-1/2} N_1^{-2}  \|v_N\|_{L^\infty_T L^2_x}\|u_{N_1}\|_{X_T^{0,1}}\| u_{N_1}\|_{L_T^\infty L^2_x} \Big)\\
&\lesssim  \sum_{N \gg 1}\sum_{N_1 \gtrsim N}\Big( \langle N_1 \rangle^{2(s-1)} N_1^{-\frac{1}{8}}\|v_N\|_{L^\infty_TL^2_x}  \| u_{N_1}\|_{L^\infty_TL^2_x}^2\\
& \quad  + {N^{3/2}}\langle N_1 \rangle^{2(s-1)}N_1^{-1} \|{v}_N\|_{X_T^{-1,1}}\| u_{N_1}\|_{L_T^\infty L^2_x}^2 \\
&\quad +N^{-1/2} \langle N_1 \rangle^{2(s-1)} N_1^{-2}  \|v_N\|_{L^\infty_T L^2_x}\|u_{N_1}\|_{X_T^{-1,1}}\| u_{N_1}\|_{L_T^\infty L^2_x} \Big)\\
&\lesssim  \sum_{N \gg 1}\sum_{N_1 \gtrsim N}\Big(  N_1^{-\frac{1}{8}}\|v_N\|_{L^\infty_TL^2_x}  \| u_{N_1}\|_{L^\infty_TH^{s-1}_x}^2\\
& \quad  + {N^{3/2}}N_1^{-1} \|{v}_N\|_{X_T^{-1,1}}\| u_{N_1}\|_{L_T^\infty H^{s-1}_x}^2 \\
&\quad +N^{-1/2} \langle N_1 \rangle^{-1} N_1^{-2}  \|v_N\|_{L^\infty_T L^2_x}\|u_{N_1}\|_{X_T^{s-1,1}}\| u_{N_1}\|_{L_T^\infty H^{s-1}_x} \Big)  \\
& \lesssim \|v\|_{\M_T^{\frac{1}{2}+}} \|u\|_{\M_T^{s-1}}\| u\|_{L_T^\infty Z^{s-1}_x}.
\end{align*}
Similarly
\begin{align*}
I^{(2)}_{N \gg 1} &= \sum_{N \gg 1}\sum_{N_1 \gtrsim N} N_1 \langle N_1 \rangle^{2(s-1)} \sup_{t \in (0,T)} |I_t(v_{\sim N_1},u_{N}, u_{N_1})| \\
&\lesssim  \sum_{N \gg 1}\sum_{N_1 \gtrsim N} N \langle N_1 \rangle^{2(s-1)} \Big(  N^{-1} N_1^{-\frac{1}{8}}\|v_{N_1}\|_{L^\infty_TL^2_x}  \| u_{N}\|_{L^\infty_TL^2_x}\|u_{N_1}\|_{L^\infty_TL^2_x}\\
& \quad  + {N^{-1/2}}N_1^{-2} \|{v}_N\|_{X_T^{0,1}} \| u_{N}\|_{L^\infty_TL^2_x}\| u_{N_1}\|_{L_T^\infty L^2_x} \\
&\quad + N^{-1/2} N_1^{-2}  \|v_N\|_{L^\infty_T L^2_x}\left(\|u_{N_1}\|_{X_T^{0,1}}\| u_{N_1}\|_{L_T^2 L^2_x} + \|u_{N_1}\|_{X_T^{0,1}}\| u_{N_1}\|_{L_T^2 L^2_x}\right) \Big)\\
& \lesssim \|v\|_{\M_T^{s}} \|u\|_{\M_T^{s-1}}\| u\|_{L_T^\infty Z^{-\frac{1}{2}+}_x} + \|v\|_{\M_T^{s}} \|u\|_{\M_T^{-\frac{1}{2}+}}\| u\|_{L_T^\infty Z^{s-1}_x}.
\end{align*}
Gathering the estimates above we have

\begin{align*}
I &  \lesssim I_1 + I_2 \\
& =   (I^{(1)}_{N \lesssim 1}+  I^{(2)}_{N \lesssim 1}) + (I^{(1)}_{N \gg 1} +   I^{(2)}_{N \gg 1}) \\
& \lesssim \|v\|_{L^\infty_T L^2_x}\|u\|_{L^2_TZ^{s-1}_x}^2 + \|v\|_{L^\infty_T Z^s_x}\|u\|_{L^\infty_TZ^{-\frac{1}{2}}_x}\|u\|_{L^\infty_TZ^{s-1}_x} + \|v\|_{\M_T^{\frac{1}{2}+}} \|u\|_{\M_T^{s-1}}\| u\|_{L_T^\infty Z^{s-1}_x}\\
& \quad + \|v\|_{\M^{s}} \|u\|_{\M_T^{s-1}}\| u\|_{L_T^\infty Z^{-\frac{1}{2}+}_x} + \|v\|_{\M_T^{s}} \|u\|_{\M_T^{-\frac{1}{2}+}}\| u\|_{L_T^\infty Z^{s-1}_x} \\
&=\left(\|v\|_{L^\infty_T Z^s_x}\|u\|_{L^\infty_TZ^{-\frac{1}{2}}_x} + \|v\|_{\M_T^{\frac{1}{2}+}} \|u\|_{\M_T^{s-1}} +\|v\|_{\M^{s}} \|u\|_{\M_T^{-\frac{1}{2}+}}+ \|v\|_{L^\infty_T L^2_x}\|u\|_{L^2_TZ^{s-1}_x}\right)\|u\|_{L^2_TZ^{s-1}_x} \\
&\qquad  + \|v\|_{\M_T^{s}} \|u\|_{\M_T^{s-1}}\| u\|_{L_T^\infty Z^{-\frac{1}{2}+}_x} \\
& \lesssim \left(\|v\|_{\M_T^{\frac{1}{2}+}} \|u\|_{\M_T^{s-1}} +\|v\|_{\M_T^{s}} \|u\|_{\M_T^{-\frac{1}{2}+}}\right)\|u\|_{L^2_TZ^{s-1}_x} + \|v\|_{\M_T^{s}} \|u\|_{\M_T^{s-1}}\| u\|_{L_T^\infty Z^{-\frac{1}{2}+}_x}.
\end{align*}
Recalling that $u= u_1 - u_2$ and $v= u_1 + u_2$ we have
\begin{align}\label{estimativa do I para dif}
\nonumber I\lesssim & \left(\|u_1 + u_2\|_{\M_T^{\frac{1}{2}+}} \|u_1 - u_2\|_{\M_T^{s-1}} +\|u_1 + u_2\|_{\M_T^{s}} \|u_1 - u_2\|_{\M_T^{-\frac{1}{2}+}}\right)\|u_1 - u_2\|_{L^2_TZ^{s-1}_x}\\
\nonumber& + \|u_1 + u_2\|_{\M_T^{s}} \|u_1 - u_2\|_{\M_T^{s-1}}\| u_1 - u_2\|_{L_T^\infty Z^{-\frac{1}{2}+}_x}\\
& \lesssim \|u_1 + u_2\|_{\M_T^{s}} \|u_1 - u_2\|_{\M_T^{s-1}}^2,
\end{align}
since $\M_T^{s}\hookrightarrow \M_T^{\frac{1}{2}+}$ and $\M_T^{s-1}\hookrightarrow \M_T^{-\frac{1}{2}+}$ because $s>\frac{1}{2}$.  Next, we estimate $J$. Following a similar approach as in the proof of Theorem \ref{des de energia}, combined with Lemma \ref{lemma do sobolev}, we obtain
\begin{align}\label{estimativa do J para dif}
\nonumber J &\lesssim  \|uv\|_{L^\infty_TH^{s-1}_x} \|\partial_x^{-1}u\|_{L^\infty_T L^2_x} \\
\nonumber &\lesssim \|u\|_{L^\infty_TH^{s-1}_x}  \|v\|_{L^\infty_TH^{s}_x}\|\partial_x^{-1}u\|_{L^\infty_T L^2_x}\\
\nonumber &\lesssim \|v\|_{\M_T^{s}} \|u\|_{\M_T^{s-1}}^2\\
&= \|u_1 + u_2\|_{\M_T^{s}} \|u_1 - u_2\|_{\M_T^{s-1}}^2.
\end{align}
Finally, summing over $N\in 2^\Z$ in inequality \eqref{desigualdade u-v}, using \eqref{estimativa do I para dif}, \eqref{estimativa do J para dif} and  applying the Theorem \ref{Littlewood-Paley Theorem} we proved that
$$
\|u_1-u_2\|^s_{L^\infty_T Z^{s-1}_x} \lesssim \|u_1^0-u_2^0\|^2_{Z^{s-1}_x} + \|u_1 + u_2\|_{\M_T^{s}}\|u_1 - u_2\|^2_{\M_T^{s-1}}.
$$
\qed
\subsection{Local Theory}\label{existecia de soluções}
This subsection is dedicated to the proof of the Theorem \ref{Teorema de existencia}. Firstly, to overcome the lack of scaling in the equation \eqref{BO2}, we will use an argument introduced by Zaiter in \cite{MR2333658} for the Ostrovsky equation. Therefore, suppose $u \in L^\infty_TZ^s_x$ be a solution of IVP \eqref{BO2} and consider $\textbf{u}(t,x):= \lambda u(\lambda^2t, \lambda x)$ for  $0< \lambda \ll 1$. Note that $\textbf{u} \in L_{\lambda^2T}^\infty Z^s_x$ and satisfies the IVP
\begin{equation}\label{BO5}
	\left\{
	\begin{array}{l}
		\partial_t \textbf{u} + \mathcal{H}\partial_x^2 \textbf{u} + \textbf{u}\partial_x \textbf{u}=  \lambda^4 \gamma \partial_x^{-1} \textbf{u}, \quad (t,x)  \in \R \times \R\\
		\textbf{u}(0,x)=\textbf{u}_0(x):=\lambda u_0(\lambda x).
	\end{array}
	\right.
\end{equation}
Given that $\|\textbf{u}_0\|_{Z^s_x} = \|\textbf{u}_0\|_{H^s_x} + \|\partial^{-1}_x \textbf{u}_0\|_{L^2_x} = \lambda(1 + \lambda^s)\|u_0\|_{Z^s_x}$ and that $\lambda > 0$ remains fixed on the right-hand side of the equation governing the evolution of problem \eqref{BO5}, we can simplify the analysis by reducing the problem to the case of small initial data. That is, we may assume that
$$
\|\textbf{u}_0\|_{Z^s_x} \ll 1.
$$ 
To simplify the notation, we will write $ \textbf{u}_0 \equiv u_0 $ and $ \textbf{u} \equiv u $ from now on. Let $ u_0 \in Z^s_x(\R) $ with $ s > \frac{1}{2} $ such that $ \|u_0\|_{Z^s_x} \ll 1 $. For $ N \in 2^\mathbb{Z} $, define $ u_0^N = P_{\leq N} u_0 \in Z^\infty_x(\R)$, and a simple computation shows that
$$
\|P_{\leq N} u_0 \|_{Z^s_x} \leq \|u_0 \|_{Z^s_x} \ll 1.
$$
We consider the problem
\begin{equation} \label{B04}
	\partial_t u + \mathcal{H} \partial_x^2 u + u\partial_xu=  \lambda^4 \partial_x^{-1} u \quad (t,x)  \in \R \times \R,
\end{equation}
 with $u(0,x)= u_0^N(x)$. It is well known from the local theory in \cite{MR2130214} (see also \cite{MR1044731}) that for $ u_0^N \in Z^\infty_x(\R) $, there exists $ T' > 0 $ and a unique solution $ u^N \in C_{T'} Z^\infty_x $ to the IVP \eqref{B04}. Fixing $s> \frac{1}{2}$, then by Theorem \ref{des de energia} we have
	\begin{align}\label{energia da aprox}
	   \nonumber \|u^N\|^2_{L_T^\infty Z^s_x} &\lesssim \|u_0^N\|^2_{Z^s_x} + \|u^N\|_{L^\infty_T L^2_x} \|u^N\|_{L^\infty_T Z^s_x}^2 \\
     &\quad + \left(1+\|u^N\|_{L^\infty_T Z^{\frac{1}{2}+}_x}\right)^2 \|u^N\|_{L^\infty_T Z^{\frac{1}{2}+}_x} \|u^N\|_{L^\infty_T Z^s_x}^2 + \|u^N\|_{L^\infty_T Z^s_x}^3,
	\end{align}
for any $T \in \left(0, \min\{1,T^\prime\}\right)$. Denoting $T^{\max}$ the maximal time of existence of function $u^N$ and by continuity of the function $T \mapsto  \|u^N\|_{L_T^\infty Z^{\frac{1}{2}+}_x}$ on $(0,T^{\max})$ we infer that
$$
\|u^N\|_{L^\infty_{T^{\prime \prime}}Z^{\frac{1}{2}+}_x} \lesssim \|u^N_0\|_{Z^{\frac{1}{2}+}_x} \ll 1,
$$ 
with $ T^{\prime \prime} := \min\{1, T^{\max}\} $. By inequality \eqref{energia da aprox} it follows that
\begin{equation}\label{desigualdade de uN}
\|u^N\|_{L^\infty_{T^{\prime\prime}}Z^{s}_x} \lesssim \| u_0^N\|_{Z^s_x} \ll 1,
\end{equation}
 for $s>1/2$. It ensures that $ T^{\prime \prime} < T^{\max} $. Consequently, we conclude that $ T^{\prime \prime} = 1 $, which implies  
\begin{equation}\label{tempo existencia}
T^{\max} \geq 1.
\end{equation}
Note that inequality \eqref{desigualdade de uN} ensures that the sequence $(u^N)_{N \in 2^\Z}$ is a bounded sequence in $L^\infty_{T}Z^{s}_x$. By a compactness argument, we obtain the existence of a solution $u \in L^\infty_{T}Z^s_x$ to the equation \eqref{B04} with $u(0)= u_0$ such that
$$
u^{N} \rightharpoonup u \quad \text{in}\,\, L^\infty_T Z^s_x.
$$
Now, let $u^1, u^2 \in L^\infty_TZ^s_x$ be two solutions of \eqref{B04} with the same initial data $u_0 \in Z^s_x(\R)$. By Theorem \ref{diferença em Zs}, Lemma \ref{lemma 3.1} and inequality \eqref{desigualdade de uN} we obtain
$$
\|u^1 - u^2 \|_{L^\infty_TZ^{s-1}_x} \lesssim \|u^1 - u^2 \|_{\M_T^{s-1}} \lesssim \|u_0 - u_0 \|_{Z^{s-1}_x} =0.
$$
Therefore, $u^1 \equiv u^2 \equiv u$. This ensures that the solution obtained is unique in the class $L^\infty_T Z^s_x \hookrightarrow L^\infty_T Z^{s-1}_x$. 

It is important to note that, based on previous estimates, the existence of a solution $ u \in L^\infty_T Z^s_x $ associated with small initial data and maximal time of existence $ T^{\max} > 1 $ guarantees the existence of solutions for arbitrarily large initial data, with a maximal time of existence satisfying:
$$
T^{\max} > (1 + \|u_0\|_{Z^s_x})^{-2}.
$$

Now, for the continuity of solution $u$ and the continuity of the flow map on $Z^s_x(\R)$ we use the standard argument given by Bona-Smith in \cite{MR385355}. For this consider $1\leq N_1 \leq N_2$ dyadic numbers and define
$$
v= u^{N_1} - u^{N_2},
$$
where $u^{N_i}$ is a solution of equation \eqref{B04} with $u^{N_i}(0)= u_0^{N_i}$ for $i=1,2$. Observe that $v$ is a solution of the equation
\begin{equation}\label{equação de u1-u2}
\partial_t v + \mathcal{H}\partial^2_x v - \frac{1}{2} \partial_x (v^2) - \partial_x (u^{N_1}v) = \partial_x^{-1}v,
\end{equation}
with $v(0)=u_0^{N_1} - u_0^{N_2}$. As a consequence of previous estimates, we easily obtain the following inequality
\begin{align}\label{3.37}
\nonumber \|v\|^2_{L_T^\infty Z^s_x} \lesssim \|u_0^{N_1}-u_0^{N_2}\|^2_{Z^s_x} &+\|v\|_{\M_T^{s}}^3+\|u^{N_1}\|_{\M_T^{s}}\|v\|_{\M_T^{s}}^2 \\
&+\|u^{N_1}\|_{\M_T^{s+1}}\|v\|_{\M_T^{s}}\|v\|_{\M_T^{s-1}}.
\end{align}
since $s>\frac{1}{2}$. Now, observe that for $s_1 \geq 0$ and $i=1,2$ we have
\begin{align*}
\|u^{N_i}(t)\|_{H^{s + s_1}_x} ^2&= \int_\R \langle \xi \rangle^{2(s + s_1)} |\F_x (u^{N_i})(t, \xi)|^2 \, d\xi \\
&=\int_\R \langle \xi \rangle^{2s}\langle \xi \rangle^{2s_1} |\eta(N_i^{-1} \xi)|^2|\widehat{u}(t, \xi)|^2 \, d\xi \\
& \lesssim\langle N_i \rangle^{2s_1} \int_\R \langle \xi \rangle^{2s}|\widehat{u}(t, \xi)|^2 \, d\xi  \\
&\lesssim \langle N_i \rangle^{2s_1}  \|u(t)\|_{H^s_x}^2\\
%& \lesssim  \langle N_i \rangle^{2s_1} \|u(t)\|_{Z^s_x}^2\\
&\lesssim \langle N_i \rangle^{2s_1}  \|u_0\|_{Z^s_x}^2.
\end{align*}
Thus $\|u^{N_i}(t)\|_{H^{s + s_1}_x} \lesssim \langle N_i \rangle^{s_1}  \|u_0\|_{Z^s_x}$. Similarly,
$$
\|(\partial_x^{-1}u^{N_i})(t)\|_{H^{s + s_1}_x} \lesssim\langle N_i \rangle^{s_1} \|u_0\|_{Z^s_x}.
$$
Therefore 
\begin{equation*}
\|u^{N_i}\|_{L^{\infty}_TZ^{s + s_1}_x} \lesssim \langle N_i \rangle^{s_1} \|u_0\|_{Z^s_x},
\end{equation*}
for all $s_1 \geq 0$ . Since $s> \frac{1}{2}$ it yields
\begin{equation}\label{3.35}
\|u^{N_i}\|_{\M^{s +s_1}_T} \lesssim \langle N_i \rangle^{s_1} \|u_0\|_{Z^s_x},
\end{equation}
for all $s_1 \geq 0$ and $i=1,2$. By Lemma \ref{lemma 3.1} and inequalities \eqref{3.35} and \eqref{3.37} we get

\begin{align*}
\|v\|^{2}_{\M^s_T} & \lesssim \left( 1+ \|u^{N_1} + u^{N_2}\|_{L^\infty_TZ^s_x}\right) \|v\|_{L^\infty_T Z^s_x} \\
& \lesssim \left( 1+ \|u^{N_1} + u^{N_2}\|_{L^\infty_TZ^s_x}\right) \Big(\|u_0^{N_1}-u_0^{N_2}\|^2_{Z^s_x} +\|v\|_{\M_T^{s}}^3  +\|u^{N_1}\|_{\M_T^{s}}\|v\|_{\M_T^{s}}^2 \\
&\hspace{8,10cm} +\|u^{N_1}\|_{\M_T^{s+1}}\|v\|_{\M_T^{s}}\|v\|_{\M_T^{s-1}} \Big)\\
&\lesssim \left(1 + \|u_0\|^2_{Z^s_x}\right) \Big(\|u_0^{N_1}-u_0^{N_2}\|^2_{Z^s_x}+\|u_0\|_{Z^s_x}\|v\|_{\M_T^{s}}^2\\
& \hspace{6,65cm} +\langle N_1 \rangle\|u_0\|_{Z^{s}_x}\|v\|_{\M_T^{s}}\|v\|_{\M_T^{s-1}} \Big) \\
& \lesssim \|u_0^{N_1}-u_0^{N_2}\|^2_{Z^s_x}+\|v\|_{\M_T^{s}}^2 +\langle N_1 \rangle\|v\|_{\M_T^{s}}\|v\|_{\M_T^{s-1}} \\
& \lesssim \|u_0^{N_1}-u_0^{N_2}\|^2_{Z^s_x}+\langle N_1 \rangle^2\|v\|_{\M_T^{s-1}}^2.
\end{align*}
But from Lemma \ref{lemma 3.1} and Theorem  \ref{diferença em Zs} we have
\begin{align*}
\|v\|_{\M_T^{s-1}}^2 & = \|u^{N_1} - u^{N_2}\|_{\M_T^{s-1}}^2\\
& \lesssim \left(1 +\|u^{N_1}+u^{N_2}\|_{L^\infty_T Z^s_x}\right)^2\|u^{N_1}- u^{N_2}\|_{L^\infty_T Z^{s-1}_x}^2 \\
& \lesssim  \|u^{N_1}_ 0-u^{N_2}_0\|^2_{Z^{s-1}_x}\\
&=  \|P_{\leq N_1}u_ 0-P_{\leq N_2}u_0\|^2_{Z^{s-1}_x} \\
&\lesssim \|P_{\leq N_1}u_ 0-u_0\|^2_{Z^{s-1}_x}.
\end{align*}
Using the definition of $P_{\leq N_1},$
\begin{align*}
\|P_{\leq N_1}u_ 0-u_0\|^2_{Z^{s-1}_x} & = \int_\R \langle \xi \rangle^{2(s -1)} |\F_x (P_{\leq N_1}u_ 0)(\xi) - \widehat{u_0}(\xi)|^2 \, d\xi\\
& \quad + \int_\R \langle \xi \rangle^{2(s -1)}|\xi|^{-2} |\F_x (P_{\leq N_1}u_ 0)(\xi) - \widehat{u_0}(\xi)|^2 \, d\xi\ \\
&= \int_\R \langle \xi \rangle^{2s} \langle \xi \rangle^{-2}|(\eta(N_1^{-1}\xi)-1)\widehat{u_0}(\xi)|^2 \, d\xi \\
&\quad + \int_\R \langle \xi \rangle^{2s} \langle \xi \rangle^{-2}||\xi|^{-2}(\eta(N_1^{-1}\xi)-1)\widehat{u_0}(\xi)|^2 \, d\xi\\
&\lesssim \langle N_1\rangle^{-2} \left(  \int_\R \langle \xi \rangle^{2s} |\F_x (P_{>N_1}u_ 0)(\xi)|^2 \, d\xi + \int_\R \langle \xi \rangle^{2s}|\xi|^{-2} |\F_x (P_{>N_1}u_ 0)(\xi)|^2 \, d\xi \right)\\
&= \langle N_1\rangle^{-2}\|P_{> N_1} u_0 \|^2_{Z^s_x}.
\end{align*}
Therefore
\begin{align}\label{limite em N_1}
\nonumber \|v\|^{2}_{\M^s_T}&\lesssim \|u_0^{N_1}-u_0^{N_2}\|^2_{Z^s_x}+\langle N_1 \rangle^2\|v\|_{\M_T^{s-1}}^2 \\
& \lesssim \|P_{> N_1} u_0 \|^2_{Z^s_x}.
\end{align}
Taking the limit as $N_1 \to 0$ in \eqref{limite em N_1} yields $\|v\|_{\M^s_T} = \|u^{N_1} - u^{N_2}\|_{\M^s_T} \to 0$, meaning that $(u^N)_{N \in 2^\Z}$ is a Cauchy sequence in the space $C_T Z^s_x$. Therefore, it converges in the $C_T Z^s_x$-norm to a solution $\bar{u}$ of \eqref{B04} with $\bar{u}(0) = u_0$. By uniqueness, we conclude that $\bar{u} \equiv u $ in $C_T Z^s_x$.
\qed
\subsection{Continuity of the Data-Solution Map}\label{continuidade do mapa dado-soluçao}

Now consider $\Phi: Z^s_x \rightarrow C_TZ^s_x$ the data-solution map. Consider $(u^n_0)_{n \in \N}$ a sequence in $Z^s_x(\R)$ such that
$$
u_0^n \rightarrow u_0 \quad \text{in} \,\, Z^s_x(\R).
$$
Let $\Phi (u^n_0):= u^n \in C_TZ^s_x$ be the solution of the IVP \eqref{BO2}. Note that, for $N \in 2^{\Z}$
\begin{equation}\label{des triangular}
\|u - u^n\|_{L_T^\infty Z^s_x} \leq \|u - P_N u\|_{L_T^\infty Z^s_x} + \|P_Nu - P_N u^n\|_{L_T^\infty Z^s_x} + \|P_Nu^n - u^n\|_{L_T^\infty Z^s_x} 
\end{equation}
where $\Phi(u_0):= u \in C_TZ^s_x$ the solution of \eqref{BO2} with $u(0) = u_0$. By the previous argument, we obtain
\begin{align*}
\|u - P_N u\|_{L_T^\infty Z^s_x} + \|u^n - P_N u^n\|_{L_T^\infty Z^s_x} & \lesssim \|u - P_N u\|_{\mathcal{M}^s_T} + \|u^n - P_N u^n\|_{\mathcal{M}^s_T}\\
&\lesssim \|P_{>N}u_0\|_{Z^s_x} + \|P_{>N}u_0^n\|_{Z^s_x}.
\end{align*}
Furthermore, following the same estimates made previously
\begin{align*}
\|P_N u - P_N u^n\|_{L_T^\infty Z^s_x}& \lesssim \|P_ Nu - P_N u^n\|_{\mathcal{M}^s_T} \\
&\lesssim  \|P_ Nu_0 - P_N u^n_0\|_{Z^s_x} + \langle N\rangle\|P_ Nu - P_N u^n\|_{\mathcal{M}^{s-1}_T}  \\
& \lesssim \|P_ Nu_0 - P_N u^n_0\|_{Z^s_x} + \langle N\rangle\|P_ Nu_0 - P_N u^n_0\|_{Z^{s-1}_x}\\
& \lesssim \|P_ Nu_0 - P_N u^n_0\|_{Z^s_x}.
\end{align*}
Substituting the above estimates in \eqref{des triangular} we obtain
\begin{align*}
\|u - u^n\|_{L_T^\infty Z^s_x} &\leq \|u - P_N u\|_{L_T^\infty Z^s_x} + \|P_Nu - P_N u^n\|_{L_T^\infty Z^s_x} + \|P_Nu^n - u^n\|_{L_T^\infty Z^s_x} \\
& \lesssim \|P_{>N}u_0\|_{Z^s_x} + \|P_{>N}u_0^n\|_{Z^s_x} + \|P_ Nu_0 - P_N u^n_0\|_{Z^s_x}\\
& \ll \varepsilon,
\end{align*}
for $N,n$ sufficiently large. Thus,
$$
\Phi(u_0^n) \rightarrow \Phi(u_0), \,\, \text{in} \,\, C_TZ^s_x ,\,\, 
$$
whenever $u^n_0 \rightarrow u_0$ in $Z^s_x(\R)$. Thus proving the continuity of the map $\Phi$.
\qed

\subsection{Global Theory}\label{global theory for MRBO}
To extend the local theory globally, we will employ certain quantities conserved by the BO-flow. First, we assume that $u$ is sufficiently smooth. However, based on the continuity properties derived in Theorem \ref{Teorema de existencia}, we observe that they indeed hold for solutions in $Z^1_x(\R)$ as long as the solution exists. Additionally, henceforth, let $C>0$ denote a constant that may change throughout the text.

We proceed as follows. Multiply the equation in \eqref{BO2} by $u$ and integrate with respect to $x$ and using the antisymmetry of the operators $\partial_x^{-1}$ and $\mathcal{H}$
\begin{equation*}
\frac{1}{2} \frac{d}{dt} \int_\R u^2 \, dx = -\int_\R u\mathcal{H}\partial^2_x u  \, dx - \int_\R u^2\partial_x u \, dx + \int_\R u\partial^{-1}_xu\, dx = 0.
\end{equation*}
Then, we obtain
\begin{equation}\label{estimativa L2}
\|u(t)\|_{L^2_x} = \|u_0\|_{L^2_x},
\end{equation}
for all $t \in \R$, where $u(0,\cdot)= u_0(\cdot)$. Next, we derive an \textit{a priori} estimate for $\|\partial_x u \|_{L^2}$  as it was done in \cite{MR1769087}. We use the fourth conserved quantity associated to the BO-flow, that is,
$$
\psi_4(u) := 2\int_\R (\partial_x u)^2\, dx + \frac{3}{2}\int_\R u^2 \mathcal{H}(\partial_x u) \, dx + \frac{1}{4}\int_\R u^4 \, dx.
$$
Now, applying the operator $F:=4\partial_x u \partial_x -u\mathcal{H}(\partial_x u)- \frac{3}{2}u^2\mathcal{H}\partial_x + u^3$ to the equation \eqref{BO2} and integrate by parts we obtain
$$
\frac{d}{dt} \psi_4(u)(t) = \int_\R (-u\mathcal{H}(\partial_x u)\partial_x^{-1}u - \frac{3}{2}u^2\mathcal{H}u + u^3\partial^{-1}_x u)\, \, dx,
$$
We need to estimate the right-hand side of the last identity, i.e.
\begin{align*}
\left|\int_\R (-u\mathcal{H}(\partial_x u)\partial_x^{-1}u - \frac{3}{2}u^2\mathcal{H}u + u^3\partial^{-1}_x u)\, \, dx \right| &\leq C \left( |I_1| + |I_2| + |I_3|\right).
\end{align*}
where
$$
I_1:=-\int_\R u\mathcal{H}(\partial_x u)\partial_x^{-1}u \, dx , \quad I_2:=- \frac{3}{2} \int_\R u^2\mathcal{H}u \, dx \,\,\, \text{and} \,\,\, I_3 := \int_\R u^3\partial^{-1}_x u\, \, dx.
$$
Let us start by estimating $I_1$. By Cauchy-Schwarz, Gagliardo-Nirenberg, Young inequalities, and by \eqref{estimativa L2} we have
\begin{align*}
|I_1| &\leq C \|\partial_x^{-1}u\|_{L^\infty_x}  \|\mathcal{H}(\partial_x u)\|_{L^2_x}   \|u_0\|_{L^2_x} \\
& \leq C \|\partial_x^{-1}u\|_{L^2_x}^{1/2}  \| u_0\|_{L^2_x}^{3/2}   \|\partial_x u\|_{L^2_x}\\
&\leq C\left( \|\partial_x^{-1}u\|_{L^2_x}  \| u_0\|_{L^2_x}^{3}  +  \|\partial_x u\|_{L^2_x}^2 \right)\\
&\leq C\left( \|\partial_x^{-1}u\|_{L^2_x}^2 + \| u_0\|_{L^2_x}^{6}  +  \|\partial_x u\|_{L^2_x}^2 \right).
\end{align*}
Similarly, for $I_2$ we have
\begin{align*}
|I_2| &\leq C \, \|u\|_{L^\infty_x}  \|\mathcal{H}u\|_{L^2_x}   \|u_0\|_{L^2_x} \\
& \leq C \, \|u_0\|_{L^2_x}^{5/2}  \|\partial_x u\|_{L^2_x}^{1/2}   \\
&\leq C \, \left(\|\partial_x u\|_{L^2_x}^{2} + \|u_0\|_{L^2_x}^{10/3}\right).
\end{align*}
Finally, we estimate $I_3$
\begin{align*}
|I_3| &\leq \|\partial_x^{-1}u\|_{L^\infty_x} \|u\|_{L^3}^3\\
& \leq C \|\partial_x^{-1}u\|_{L^2_x}^{1/2}\|u_0\|_{L^2_x}^{5/2}\|\partial_x u\|_{L^2_x}\\
&\leq C \left(\|\partial_x^{-1}u\|_{L^2_x}\|u_0\|_{L^2_x}^{5}+ \|\partial_x u\|_{L^2_x}^2 \right)\\
&\leq C \left(\|\partial_x^{-1}u\|_{L^2_x}^2 +\|u_0\|_{L^2_x}^{10}+ \|\partial_x u\|_{L^2_x}^2 \right).
\end{align*}
Combining the estimates for $I_1$, $I_2$ and $I_3$ we obtain that
\begin{align*}
\left|\int_\R (-u\mathcal{H}(\partial_x u)\partial_x^{-1}u - \frac{3}{2}u^2\mathcal{H}u + u^3\partial^{-1}_x u)\, \, dx \right| &\leq C \big[  \|\partial_x u\|_{L^2_x}^2+ \|\partial_x^{-1}u\|_{L^2_x}^2+ \|u_0\|_{L^2_x}^{6}\\
& \qquad+ \|u_0\|_{L^2_x}^{10/3} + \|u_0\|_{L^2_x}^{10} \big].
\end{align*}
Then
% Denoting
% $$
% G(t)= G(\|\partial_x^{-1}u(t)\|_{L^2_x}):=  \|\partial_x^{-1}u\|_{L^2_x}^2+ \|u_0\|_{L^2_x}^{6}+ \|u_0\|_{L^2_x}^{10/3} + \|u_0\|_{L^2_x}^{10},
% $$
% and knowing that , we have
\begin{align}\label{2.16}
\nonumber \frac{d}{dt} \psi_4(u)(t) &\leq \left|\int_\R (-u\mathcal{H}(\partial_x u)\partial_x^{-1}u - \frac{3}{2}u^2\mathcal{H}u + u^3\partial^{-1}_x u)\, \, dx \right| \\
&\leq C \left( \|\partial_x u \|_{L^2_x}^2 +  \|\partial_x^{-1}u\|_{L^2_x}^2 + C_1(\|u_0\|_{L^2_x}) \right).
\end{align}
It is not difficult to see that
$$
\psi_4(u(0)) \leq {C}_2(\|u_0\|_{Z^1_x}). 
$$
Integrating \eqref{2.16} with respect to t and using  it follows then that
\begin{equation}\label{2.17}
    \psi_4(u)(t) \leq C \left({C}_3(\|u_0\|_{Z^1_x})+  \int_0^t \|\partial_x^{-1}u(\tau)\|_{L^2_x}^2\,d\tau + \int_0^t \|\partial_x u(\tau) \|_{L^2_x}^2 \, d\tau\right).
\end{equation}
On the other hand, isolating the term $\|\partial_x u \|_{L^2_x}^2$ in $\psi_4$, we obtain
$$
\|\partial_x u \|_{L^2_x}^2 = \frac{1}{2} \psi_4(u)(t) - \frac{3}{4}\int_\R u^2 \mathcal{H}(\partial_x u) \, dx - \frac{1}{8}\int_\R u^4 \, dx.
$$
 Gagliardo-Nirenberg and Young's inequalities imply then
\begin{align}\label{2.20}
\|\partial_x u \|_{L^2_x}^2 \leq \psi_4(u)(t) + (\varepsilon_1 +\varepsilon_2)\|\partial_x u \|_{L^2_x}^2 + \tilde{C}_{\varepsilon_1,\varepsilon_2}\|u_0\|_{L^2_x}^6,
\end{align}
where $\varepsilon_1,\varepsilon_2 \in \R_+$. Combining \eqref{2.17} and \eqref{2.20}, together with suitable choices of $\varepsilon_1$ and $\varepsilon_2$ we obtain
\begin{align}\label{estimativa derivada}
\|\partial_x u \|_{L^2_x}^2 \leq C \left( C_4(\|u_0\|_{Z^1_x}) +  \int_0^t \|\partial_x^{-1}u(\tau)\|_{L^2_x}^2\,d\tau+   \int_0^t \|\partial_x u(\tau) \|_{L^2_x}^2 \, d\tau \right).
\end{align}
% By the Gronwall inequality we have
% \begin{equation}
% \|\partial_x u \|_{L^2_x}^2 \leq C  \left[C_5(\|u_0\|_{Z^1_x})+  \int_0^T\left(\int_0^t \|\partial_x^{-1}u(\tau)\|_{L^2_x}^2\,d\tau\right)\, dt\right] \exp(T),
% \end{equation}
% {for} $0 \leq t \leq T$. 
Now we estimate the quantity $\|\partial_x^{-1} u \|_{L^2_x}$. Applying the operator $\partial_x^{-1}$ in the equation \eqref{BO2} and taking the $L^2_x$-inner product with $\partial_x^{-1}u$ we obtain

$$
\frac{1}{2} \frac{d}{dt}\|\partial_x^{-1}u \|_{L^2_x}^2 + \frac{1}{2} \int_\R u^2 \partial_x^{-1}u \, dx = 0.
$$
The Cauchy-Schwarz,  Gagliardo-Nirenberg and Young's inequalities imply that
\begin{align*}
\frac{d}{dt}\|\partial_x^{-1}u \|_{L^2_x}^2 &\leq C\|u\|_{L^4_x}^2\|\partial_x^{-1}u \|_{L^2_x} \\
&\leq C \|u_0\|_{L^2_x}^{3/2}\|\partial_x u \|_{L^2_x}^{1/2}\|\partial_x^{-1}u \|_{L^2_x} \\
&\leq C \left( \|u_0\|_{L^2_x}^{3}\|\partial_x u \|_{L^2_x}+ \|\partial_x^{-1}u \|_{L^2_x}^2\right)\\
&\leq C \left( \|u_0\|_{L^2_x}^{6}+\|\partial_x u \|_{L^2_x}^2+ \|\partial_x^{-1}u \|_{L^2_x}^2\right).
\end{align*}
Therefore, we obtain the following estimate
\begin{equation}\label{estimativa antiderivada}
\|\partial_x^{-1} u \|_{L^2_x}^2\leq C \left( C_5(\|u_0\|_{Z^1_x}) +  \int_0^t \|\partial_x^{-1}u(\tau)\|_{L^2_x}^2\,d\tau+   \int_0^t \|\partial_x u(\tau) \|_{L^2_x}^2 \, d\tau \right).
\end{equation}
By \eqref{estimativa derivada}, \eqref{estimativa antiderivada}  we have
\begin{equation*}
\|\partial_x u \|_{L^2_x}^2+\|\partial_x^{-1} u \|_{L^2_x}^2\leq C \left( C_6(\|u_0\|_{Z^1_x}) +  \int_0^t \left[\|\partial_x^{-1}u(\tau)\|_{L^2_x}^2+ \|\partial_x u(\tau) \|_{L^2_x}^2 \right] \, d\tau \right).
\end{equation*}
Now, by the Gronwall inequality we can conclude that
\begin{equation}\label{estimativa anti+ derivada}
    \|\partial_x u \|_{L^2_x}^2+\|\partial_x^{-1} u \|_{L^2_x}^2\leq C_7(\|u_0\|_{Z^1_x}) \exp(CT), \quad 0\leq t \leq T.
\end{equation}

Since we have the local theory well established, the \textit{a priori} estimate \eqref{estimativa anti+ derivada} above, can be used to extend the local solution globally in $Z^1_x(\R)$, because
$$
\|u\|_{Z^1_x} ^2\sim  \left( \| u \|_{L^2_x}^2+ \|\partial_x u \|_{L^2_x}^2+\|\partial_x^{-1} u \|_{L^2_x}^2\right).
$$
As the BO equation possesses infinite conservation laws (for each order $1/2$ of regularity, for instance, in $H^1_x, H_x^{3/2}, H_x^2,...$), a similar argument to the one previously employed, along with interpolation, allows us to globally extend the solution to $Z^s_x(\mathbb{R})$ for $s \ge 1$. This completes the proof. \qed

\section{The Rotational Modified Intermediate Long Wave Equation}\label{seçao ILW}
 This section will be dedicated to proving the local well-posedness for the following Cauchy problem
\begin{equation}\label{ILW6}
	\left\{
	\begin{array}{l}
		\partial_t v + \mathcal{T}_1  \partial_x^2 v+ \partial_x v+ v\partial_xv=  \gamma \partial_x^{-1} v\quad (t,x)  \in \R \times \R,\\
		v(0,x)=v_0(x).
	\end{array}
	\right.
\end{equation}
It becomes clear to the reader that a significant portion of the previously employed techniques follows a similar framework to this model, albeit with minor adaptations. Therefore, let us recall that the symbol of the equation \eqref{ILW2} is given by
$$
\omega(\xi) = - \xi^2\coth{\xi}+ \xi- \gamma \xi^{-1}.
$$
Consequently consider now the \textit{resonance function} given by the symbol $\omega(\xi)$
\begin{equation}\label{função de ressonancia ILW}
\Lambda(\xi_1,\xi_2):=\omega(\xi_1 + \xi_2) - \omega(\xi_1) - \omega(\xi_2),
\end{equation}
for $(\xi_1,\xi_2) \in \R \times \R$. First of all, we obtain a lemma similar to Lemma \ref{controle resonant} for the resonance function \eqref{função de ressonancia ILW}.
\begin{lemma}\label{controle resonant ILW}
Consider $(\xi_1,\xi_2) \in \R^2$. If $\min \{|\xi_1|, |\xi_2|\} \gg 1$ then
$$
|\Lambda(\xi_1,\xi_2)|  \lesssim |\xi|^2_{\max}|\xi|_{\min} + |\xi|_{\max}|\xi|^{-2}_{\min}.,
$$
where $|\xi|_{\max}= \max \{|\xi_1|,|\xi_2|, |\xi_1 + \xi_2|\}$ and $|\xi|_{\min}= \min \{|\xi_1|,|\xi_2|, |\xi_1 + \xi_2|\}.$
\end{lemma}
\dem 
Again, without loss of generality, assume $\gamma=1$. From \cite{MR4215711}, the symbol $\omega$ is a real odd function such that 
\begin{equation}\label{KdV behavior}
    \omega(\xi) = \frac{1}{3}\xi^3 + O(\xi^5) \quad \mbox{as} \quad |\xi| \rightarrow 0
\end{equation}
and 
\begin{equation}
    \omega(\xi) = \xi|\xi| - \xi - \xi^{-1}+ O(1) \quad \mbox{as} \quad  |\xi| \rightarrow \infty.
\end{equation}
Therefore we can separate the proof of the result in two cases: 

\noindent \textit{\underline{Case 1}}: $|\xi_2| \gtrsim  |\xi_1|$.

By applying a straightforward method similar to that in Lemma \ref{controle resonant}, we can derive this result. The case $|\xi_1| \gtrsim  |\xi_2|$ follows by symmetry.

\noindent \textit{\underline{Case 2 }}: $|\xi_2| \sim  |\xi_1|$  .

Here is the case where the singularity appears, then we need to be carefully. By property \eqref{KdV behavior}, there exist $\varepsilon > 0 $ and $M > 0, $ such that, if $|\xi_1 + \xi_2|< \varepsilon, $ then 
\begin{equation}
    \begin{aligned}
     |\Lambda(\xi_1,\xi_2)| & \lesssim |\xi_1 + \xi_2|^3 + M\varepsilon^5 + \left|  \frac{1}{\xi_1 +\xi_2} - \frac{1}{\xi_1} - \frac{1}{\xi_2} \right| \\
     & \lesssim |\xi|^2_{\max}|\xi|_{\min} + |\xi|_{\max}|\xi|^{-2}_{\min}.
    \end{aligned}
\end{equation}
Using that  $\min \{|\xi_1|, |\xi_2|\} \gg 1$, the result follows.
\qed

From Lemma \ref{controle resonant ILW}  and Bernstein's inequalities, we can obtain a similar result obtaneid in Lemma \ref{lemma 2.2}.
\begin{lemma}\label{lemma 2.2 ILW}
	Let $N_i \in 2^\Z$, $L_i \in 2^\N$ with $i=1,2,3$ and $0<N_1 \leq N_2 \leq N_3$. Furthermore consider $f_i \in \Sc(\R^2)$ such that
$$
\text{supp }\widehat{f_i}\subset \{\xi \in \R: |\xi| \sim N_i \},
$$
for each $i=1,2,3$. Then
$$
\int_\R\int_\R (\Pi(\bar{Q}_{L_1}f_1,\bar{Q}_{L_2}f_2) \bar{Q}_{L_3}f_3)(t,x)\,\, dt\, dx =0,
$$
for all pseudoproduct $\Pi \equiv \Pi_{\Psi}$, whenever the $\displaystyle \max_{i=1,2,3}L_i \gtrsim N_1N_2^2 +N_1^{-2}N_2$  is not satisfied.\qed
\end{lemma}

\subsection{Local Theory}
Note that, as a consequence of Lemmas \ref{controle resonant ILW} and \ref{lemma 2.2 ILW}, Theorems \ref{des de energia} and \ref{diferença em Zs} also hold for the problem \eqref{ILW2}. Therefore, we can now proceed to establish the local well-posedness for this IVP. Let us suppose $v \in L^\infty_T Z^s_x$ is a solution of IVP \eqref{ILW2}, and consider $\textbf{v}(t,x):= \lambda v(\lambda^2t, \lambda x)$ for $0< \lambda \ll 1$. It is noteworthy that $\textbf{v} \in L^\infty_{\lambda^2 T} Z^s_x$ and satisfies the IVP
\begin{equation}\label{BO3}
	\left\{
	\begin{array}{l}
		\partial_t \textbf{v} + \mathcal{T}_1\partial_x^2 \textbf{v}+ \partial_x \textbf{v} +  \textbf{v}\partial_x\textbf{v}=  \lambda^4 \gamma\partial_x^{-1} \textbf{v}, \quad (t,x)  \in \R \times \R, \\
		\textbf{v}(0,x)=\textbf{v}_0(x):=\lambda v_0(\lambda x).
	\end{array}
	\right.
\end{equation}
Since $\|\textbf{v}_0\|_{Z^s_x} = \|\textbf{v}_0\|_{Z^s_x} + \|\partial^{-1}_x \textbf{v}_0\|_{L^2_x} = \lambda(1 +\lambda^s)\| v_0\|_{Z^s_x}$ and $\lambda>0$ will remain fixed on the right side of equation in the evolution of problem \eqref{BO3}, we can reduce our problem to the case of small initial data, \textit{i.e.} we can assume that
$$
\|\textbf{v}_0\|_{Z^s_x} \ll 1.
$$ 
Again, in order to simplify the notation, we will write $ \textbf{v}_0 \equiv v_0 $ and $ \textbf{v} \equiv v $ from now on.

By following the same arguments used in \cite{MR2130214} (possibly combined with the techniques developed in \cite{MR1044731}), it is possible to ensure the existence of smooth solutions for \eqref{ILW2}, assuming regular initial data, with a time of existence $T > 0$. Given that the energy estimate \eqref{des de energia} holds (following the same approach), and following the same steps outlined in Subsection \ref{existecia de soluções}, we conclude that the problem \eqref{ILW2} is locally well-posed in $Z^s_x(\R)$, provided that $s > 1/2$.

Furthermore, given that the estimate between two solutions \eqref{diferença em Zs} holds, it follows from the same arguments employed in Subsection \ref{continuidade do mapa dado-soluçao} that the data-solution map $v_0 \mapsto C_TZ^s_x$ of the equation \eqref{ILW2} is continuous.
\qed

\subsection{Global Theory}

For the global well-posedness in $Z^s_x(\R)$ for $s\geq 1$, it suffices to note firstly that if $v(t)$ is a solution to \eqref{ILW2}, then
$$
\|v(t)\|_{L^2_x}= \|v_0\|_{L^2_x},
$$
for all $0 \leq t \leq T$, and considering the following conserved quantity associated with the ILW-flow, \textit{i.e.}
\begin{align*}
\varphi_4(u)&:= \frac{1}{2} \int_\R (\partial_x u)^2 \, dx +\frac{3}{2} \int_\R u^2\mathcal{T}_1\partial_x u \, dx + \frac{3}{2} \int_\R (\mathcal{T}_1\partial_x u)^2 \, dx \\
&\quad+ \frac{9}{2} \int_\R u\mathcal{T}_1\partial_x u \, dx + \frac{3}{2} \int_\R u^3 \, dx +\frac{3}{2} \int_\R u^2\, dx +  \frac{1}{4} \int_\R u^4 \, dx.
\end{align*}
Thus, applying the operator
$$
\Tilde{F}:= \partial_x v \partial_x + 3v\mathcal{T}_1\partial_x v + \frac{3}{2}v^2\mathcal{T}_1\partial_x + 3 \mathcal{T}_1\partial_x v\mathcal{T}_1\partial_x +\frac{9}{2}\mathcal{T}_1\partial_xv +\frac{9}{2}v\mathcal{T}_1\partial_x + \frac{9}{2}v^2 + 3v + v^3
$$
to the equation \eqref{ILW6} and integrating by parts we obtain
$$
\frac{d}{dt} \varphi_4(v)(t) = \int_\R \left(3v\mathcal{T}_1\partial_xv \partial_x^{-1}v + \frac{3}{2}v^2\mathcal{T}_1v +\frac{9}{2}v^2 \partial_x^{-1}v + v^3\partial^{-1}_x u\right)\, \, dx.
$$
By applying Gagliado-Nirenberg and Cauchy-Schwarz inequalities, integration by parts, and using the relation $\partial_x \partial_x^{-1} = \partial_x^{-1} \partial_x = I_d$, where $I_d$ denotes the identity operator, we obtain
\begin{align*}
\left| \int_\R v^2\mathcal{T}_1v \, dx\right| &\leq C \|v\|_{L^\infty_x}\|\partial_x^{-1}v\|_{L^2_x}\|\mathcal{T}_1\partial_x^{}v\|_{L^2_x} \\
&\leq C \|v_0\|_{L^2_x}^{3/2}\|\partial_x v\|_{L^2_x}^{1/2}\|\partial_x^{-1}v\|_{L^2_x}\\
&\leq C \left(\|v_0\|_{L^2_x}^{3}\|\partial_x v\|_{L^2_x}+\|\partial_x^{-1}v\|_{L^2_x}^2\right)\\
&\leq C\left(\|v_0\|_{L^2_x}^{6}+\|\partial_x v\|_{L^2_x}^2+\|\partial_x^{-1}v\|_{L^2_x}^2\right).
\end{align*}
Similarly
$$
\left|\int_\R 3v\mathcal{T}_1\partial_xv \partial_x^{-1}v \, dx \right| \leq   C\left(\|v_0\|_{L^2_x}^{6}+\|\partial_x v\|_{L^2_x}^2+\|\partial_x^{-1}v\|_{L^2_x}^2\right).
$$
With these estimates and following steps similar to those in Subsection \ref{global theory for MRBO}, we can obtain \textit{a priori} estimates for solutions of \eqref{ILW2} in $Z^1_x(\mathbb{R})$ as follows
\begin{equation*}
   \sup_{0 \leq t \leq T} \left(\|\partial_x v\|_{L^2_x}^2+\|\partial_x^{-1} v \|_{L^2_x}^2\right)\leq \tilde{C}(\|v_0\|_{Z^1_x}) \exp(CT).
\end{equation*}
These estimate can be used to extend the local solutions globally in $Z^s_x(\R)$ for $s \geq 1$ given that the ILW equation is a completely integrable system, \textit{i.e.}, it possesses infinitely many conserved quantities.
\qed

\noindent\textbf{Acknowledgments:} The authors would like to thanks Prof. Dr. Felipe Linares for all the discussions and insights about this problem.

\bibliographystyle{amsalpha}

\end{document}